\documentclass[11pt,reqno,a4paper]{amsart}

\usepackage{amsmath,amssymb,amsthm,amsaddr}
\usepackage[utf8]{inputenc}
\usepackage[british]{babel}

\usepackage{enumitem}

\usepackage{tikz}
\usepackage[framemethod=TikZ]{mdframed}

\usepackage{multirow}

\usepackage[margin=2.5cm,top=4cm,bottom=3cm,footskip=1cm,headsep=2cm]{geometry}

\usepackage{accents}
\newcommand{\ubar}[1]{\underaccent{\bar}{#1}}

\newtheorem{theorem}{Theorem}
\renewcommand{\thetheorem}{\arabic{theorem}}

\theoremstyle{plain}
\newtheorem{lemma}[theorem]{Lemma}

\theoremstyle{definition}

\newtheorem{problem}[theorem]{Problem}

\newtheorem{remark}[theorem]{Remark}

\def\eps{\varepsilon}

\def\dx{\partial_x}

\def\dt{\partial_t}

\def\V{\mathcal{V}}

\def\E{\mathcal{E}}
\def\T{\mathcal{T}}

\def\R{\mathbb{R}}

\def\AS{\mathcal{AS}}

\def\H{\mathcal{H}}
\def\D{\mathcal{D}}
\def\T{\mathcal{T}}

\def\R{\mathbb{R}}
\def\dtau{\bar d_{\tau}}
\def\ddtau{\frac{d}{d\tau}}
\def\ptau{\partial_{\tau}}
\def\pptau{\partial_{\tau\tau}}
\def\resq{\operatorname{res}_1}
\def\resr{\operatorname{res}_2}
\def\Dtau{\Delta \tau}

\def\la{\langle}
\def\ra{\rangle}

\newcommand{\hn}[1]{{#1}_h^n}
\newcommand{\hatn}[1]{\hat{#1}_h^n}

\newcommand{\LL}[1]{\|#1\|_{L^2}^2}

\usepackage[pdfa,bookmarksnumbered=true]{hyperref}
\hypersetup{
  colorlinks=true,
  allcolors = red,
  linkcolor=blue, 
  menucolor=blue,
  urlcolor = blue,
  frenchlinks=false,
  pdfborder={0 0 0},
  naturalnames=false,
  hypertexnames=false,
  breaklinks
}

\title[Asymptotic preserving discretization for gas transport]{An asymptotic-preserving discretization scheme\\ for gas transport in pipe networks}
\author{H. Egger$^*$, J. Giesselmann$^\dag$, T. Kunkel$^\dag$, \and N. Philippi$^\dag$}
\address{$^*$Johann Radon Institute for Computational and Applied Mathematics and Institute for Computational Mathematics, Johannes-Kepler University Linz, Austria}
\address{$^\dag$Department of Mathematics, TU Darmstadt, Germany}
\email{egger@numa.uni-linz.ac.at}
\email{giesselmann@mathematik.tu-darmstadt.de}
\email{tkunkel@mathematik.tu-darmstadt.de}
\email{philippi@mathematik.tu-darmstadt.de}

\begin{document}

\begin{abstract}
We consider the simulation of barotropic flow of gas in long pipes and pipe networks. Based on a Hamiltonian reformulation of the governing system, a fully discrete approximation scheme is proposed using mixed finite elements in space and an implicit Euler method in time. Assuming the existence of a smooth subsonic solution bounded away from vacuum, a full convergence analysis is presented based on relative energy estimates. 
Particular attention is paid to establishing error bounds that are uniform in the friction parameter. As a consequence, the method and results also cover the parabolic problem arising in the asymptotic large friction limit. 
The error estimates are derived in detail for a single pipe, but using appropriate coupling conditions and the particular structure of the problem and its discretization, the main results directly generalize to pipe networks. 
Numerical tests are presented for illustration. 
\end{abstract}

\maketitle

\vspace*{-1em}

\begin{quote}
 
\end{quote}

\begin{quote}
\noindent
{\small {\bf AMS-classification (2000):}
35L65, 65M15, 65M60, 76M10}
\end{quote}

\section{Introduction}

We consider the systematic numerical approximation of gas transport in long pipes and pipeline networks. The flow of gas on each pipe is described by the barotropic Euler equations with a source term accounting for friction at the pipe walls. Under standard operating conditions, the gas flows at moderate velocities and we therefore consider the following rescaled equations~\cite{Brouwer11} that account for the low Mach or high friction regime:
\begin{align}
	a \partial_\tau \rho + \dx m &= 0,\label{eq:cgas_1}\\
	\eps^2 \partial_\tau m + \dx \Big(\frac{\eps^2  m^2}{a \rho} + a p(\rho) \Big) &= -\frac{\gamma}{a \rho}|m|m.\label{eq:cgas_2}
\end{align}
Here $\rho$ is the density of the gas, $m$ is the mass flux, $a$ is the constant cross-sectional area of the pipe, $\gamma$ is a friction coefficient, $p=p(\rho)$ is the pressure, and $\eps$ is a scaling parameter proportional to the Mach number; see \cite{Brouwer11} or Appendix~\ref{app:rescale} for details.

Together with appropriate boundary conditions and, in the case of networks, with coupling conditions that ensure conservation of mass and energy at pipe junctions \cite{Egger18,Reigstad15}, the one-dimensional Euler equations \eqref{eq:cgas_1}--\eqref{eq:cgas_2} define a port-Hamiltonian system.  
The particular problem structure becomes apparent by an appropriate reformulation of the equations and a variational characterization of its solutions, which immediately leads to an energy dissipation inequality;  see \cite{EggerGiesselmann, GiesselmannLattanzioTzavaras_17} and Section~\ref{sec:single} fro details. 
In the weak formulation, the coupling and boundary conditions are incorporated variationally, which allows for a structure-preserving discretization by Galerkin projection. 
In this paper, we consider the discretization by a mixed finite element scheme in space, using piecewise constant approximations $\rho_h$ for $\rho$ and continuous and piecewise linear approximations $m_h$ for $m$, combined with an implicit Euler time discretization. 
This amounts to a standard approximation for related linear wave propagation problems \cite{Geveci88,Joly03}. A closely related scheme has been considered in~\cite{Egger18,Egger19} for isentropic flow in pipe networks; let us also refer to \cite{Cardoso19,Liljegren20} for similar approaches. 

The main contribution of the current work is a rigorous convergence analysis for the discretization scheme outlined above. We will establish order optimal error estimates 
\begin{align*} 
	\|\rho(\tau^n)-\rho_h^n\|_{L^2}^2 + \eps^2\|m(\tau^n)-m_h^n\|_{L^2}^2 +  \sum_{k=1}^n  \Dtau \|m(\tau^k)-m_h^k\|_{L^3}^3
	\leq C \big(\Dtau^2 + h^2\big) 
\end{align*}
in the temporal and spatial mesh size, which hold uniformly for all $0 \le \eps \le \bar \eps$. 
In particular, the proposed method yields a viable discretization scheme with order optimal convergence also for the parabolic limit problem, which results from formally setting $\eps =0$ in the equations.
This parabolic model  is widely used in the gas network community \cite{Brouwer11,Burlacu19}, and studied intensively in the literature, see e.g. \cite{Bamberger79} or \cite{SchoebelKroehn20} for a full analysis on networks.
The parameter-robust error estimate above is proven under the assumptions that 
the flow is subsonic and bounded away from vacuum and that 
certain derivatives of the solution are bounded uniformly in $\eps$. These assumptions are reasonable for gas flows in pipe networks under standard operating conditions.
The consideration of shocks or discontinuities, as in \cite{hertyEMS}, is beyond the scope of this work.

Our analysis is based on discrete stability of the proposed scheme, which is established via \emph{relative energy estimates} and uses similar arguments as in \cite{EggerGiesselmann}, where the asymptotic limit $\eps \searrow 0$ of \eqref{eq:cgas_1}--\eqref{eq:cgas_2} was investigated on the continuous level.
Relative energy or entropy estimates are a well known tool for the analysis of quasi-linear partial differential equations; see \cite{Jungel16} for an overview on its use in parabolic equations and \cite{Dafermos05} for applications in hyperbolic balance laws. 
They have also been employed successfully for the numerical analysis of the compressible Navier-Stokes equations \cite{Feireisl18,Gallouet16,Kwon20} and of the Euler equations in the large friction limit~\cite{Berthon_Bessemoulin_Mathis_2017}.

The remainder of this manuscript is organized as follows: 
In Section~\ref{sec:single}, we state our basic assumptions and introduce the Hamiltonian reformulation of the barotropic Euler equations as well as its variational formulation. Moreover, we briefly discuss the underlying port-Hamiltonian structure and the corresponding energy dissipation law. 
The discretization method, some basic properties, and a complete statement of the above error estimates are presented in Section~\ref{sec:disc} and their detailed proof is given in Section~\ref{sec:proof}.
Particular emphasis is payed on explicitly tracking the dependence of all estimates on the parameter~$\eps$. 
In Section~\ref{sec:network}, we  show that the proposed method and its analysis seamlessly generalize to gas networks, if appropriate coupling conditions are required at pipe junctions. 
For illustration of our results, 
some numerical tests are presented in Section~\ref{sec:numerics}. 
Some auxiliary results are summarized in a small appendix.


\section{Formulation on a single pipe} \label{sec:single}

In this section, we briefly state our basic assumptions, introduce a weak formulation of the problem,  and collect some auxiliary results concerning the energy of the system. 

\subsection{Problem statement}

We start with rewriting the momentum equation into an evolution equation for the rescaled flow velocity $w$. Then, \eqref{eq:cgas_1}--\eqref{eq:cgas_2} can be stated as 
\begin{align}
	a \partial_\tau \rho + \dx m &= 0,\label{eq:gas1}\\
	\eps^2 \partial_\tau w + \dx h &= -\gamma|w|w\label{eq:gas2}
\end{align}
for $0<x<\ell,\ \tau>0$, with mass flux $m$ and total specific enthalpy $h$ given by
\begin{align} \label{eq:gas3}
    m=a\rho w,\qquad h=\tfrac{1}{2}\eps^2w^2 +P'(\rho).
\end{align}
The potential $P$ here is connected to the pressure $p$ by $p'(\rho)=\rho P''(\rho)$; see
Appendix~\ref{app:rescale} for a detailed derivation of the equations. 
Any pair of functions
\begin{align*} 
\rho,w \in C^1([0,\tau_{max}];L^2(0,\ell)) \cap C^0([0,\tau_{max}];H^1(0,\ell)) \qquad \text{with} \qquad  \rho>0,
\end{align*}
that satisfies \eqref{eq:gas1}--\eqref{eq:gas3} in a pointwise sense, will be called a \emph{classical solution}.
For such smooth solutions away from vacuum, the two systems \eqref{eq:cgas_1}--\eqref{eq:cgas_2} and \eqref{eq:gas1}--\eqref{eq:gas3} are equivalent.
We further complement the system \eqref{eq:gas1}--\eqref{eq:gas3} by boundary conditions 
\begin{align} \label{eq:gas4}
h(v,\tau)=h^v_\partial(\tau)
\end{align}
at the end points $v\in\{0,\ell\}$ of the pipe and $\tau>0$,
which make sense for classical solutions.

\subsection{Main assumptions}

For our analysis later on, we assume that 
\begin{enumerate}[label=(A\arabic*)]
	\item \label{A1}
	the pressure potential $P:\R_+\to \R$ is smooth and strongly convex;
	\item \label{A2}
	the constants $\eps,\gamma,a$ satisfy $0\le\eps\le\bar \eps$,  $0<\ubar\gamma\le\gamma\le\bar\gamma$, and $0 < \ubar a \le a \le \bar a$, as well as
	\begin{align} \label{eq:subsonic}
	\rho P''(\rho)\ge 4 \bar\eps^2 |\bar w|^2 \qquad \forall 2\ubar \rho/3 \leq \rho \leq 3\bar \rho/2
	\end{align}
	with appropriate positive constants designated by bar symbols;
    \item \label{A3}	
    there exists a classical solution $(\rho,w)$ such that 
    \begin{align} \label{eq:bounds} 
    0 < \ubar \rho\le \rho (t,x) \le \bar \rho
    \qquad \text{and} \qquad 
    -\bar w \le w(t,x) \le \bar w
    \end{align}
     for all $0 \le \tau \le \tau_{max}$ and a.a. $0 \le x \le \ell$; this is called a \emph{subsonic bounded state solution}.
\end{enumerate}
In order to obtain quantitative error estimates, we additionally require that
\begin{enumerate}[label=(A\arabic*)]
    \setcounter{enumi}{3}
    \item \label{A4}
    the solution provided by \ref{A3} is sufficiently smooth, i.e., $\rho$ and $w$ are uniformly bounded in 
    $W^{2,\infty}(0,\tau_{max};L^2(0,\ell))$ and $W^{1,\infty}(0,\tau_{max};H^1(0,\ell))$
    by a constant $\bar C$.
\end{enumerate}
Bounds for derivatives of $m$ and $h$ follow from equation \eqref{eq:gas3} and the previous assumptions. 
\begin{remark}
Condition \ref{A2} implies strict convexity of the pressure potential $P(\rho)$, or equivalently, strict monotonicity of $p(\rho)$, which is a natural thermodynamic requirement ensuring hyperbolicity of the barotropic Euler equations.
Let us further note that $c(\rho)=\sqrt{p'(\rho)}=\sqrt{\rho P''(\rho)}$ is the speed of sound, and hence assumption \ref{A3} characterizes solutions bounded away from the sonic point and from vacuum. 
In that case, one can require exactly one boundary condition at every end of the pipe; see e.g. \cite{Yee_1981}.
\end{remark}

\subsection{Weak form of the problem}

The following variational characterization of classical solutions is the starting point for our discretization strategy and the further investigations. 
\begin{lemma} \label{lem:var}
Let $(\rho,w)$ denote a classical solution of \eqref{eq:gas1}--\eqref{eq:gas4}. Then 
\begin{align}
	\la a\partial_{\tau}\rho(\tau), q\ra + \la \dx m(\tau), q\ra &= 0 &&\forall q\in L^2(0,\ell), \label{eq:var1}\\
	\la \eps^2\partial_{\tau} w(\tau), r\ra - \la h(\tau), \dx r\ra &= - h_\partial(\tau) r \vert_0^\ell - \la \gamma |w(\tau)| w(\tau), r\ra &&\forall r\in H^1(0,\ell),  \label{eq:var2}
\end{align}
and all $0 \le \tau \le \tau_{max}$. 
Here and below, we use $\la u,v\ra := \int_0^\ell a \, b \, dx$ to denote the standard scalar product of $L^2(0,\ell)$ and $h_\partial r \vert_0^\ell := h_\partial^\ell r(\ell) - h_\partial^0 r(0)$ to abbreviate the boundary terms.
\end{lemma}
\begin{proof}
The above identities follow immediately from \eqref{eq:gas1}--\eqref{eq:gas2} by multiplying with appropriate test functions and integration over the domain $(0,\ell)$. For the second equation, we utilize integration-by-parts for the spatial derivatives and the boundary conditions \eqref{eq:gas4}. 
\end{proof}

\subsection{Energy dissipation and convexity}
A particular feature of problem \eqref{eq:gas1}--\eqref{eq:gas2} is that the quantities arising in the equations can be 
understood, respectively, as \textit{state variables} $(\rho,w)$ and \textit{co-state variables} $(h,m)$, directly linked via the associated energy functional
\begin{align} \label{eq:H}
   \H(\rho,w):=\int_{0}^{\ell} a \, \left(\tfrac{1}{2}\eps^2 \rho w^2 +P(\rho)\right)\, dx.
\end{align}
More precisely
\begin{align}\label{eq:energy-varder}
    \delta_\rho\H(\rho,w)=ah, \quad \delta_w \H(\rho,w)=\eps^2 m
\end{align}
are the \emph{variational derivatives} of $\H$, i.e., the partial derivatives of the integrand in \eqref{eq:H}.
We will write $\H' = (\delta_\rho \H, \delta_w \H)$ for the variational derivative of $\H$ in the sequel. 
As a direct consequence of these relations, we obtain the following structural property. 
\begin{lemma} \label{lem:pH}
Let $(\rho,w)$ denote a classical solution of \eqref{eq:gas1}--\eqref{eq:gas4}. Then 
\begin{align}\label{eq:energy}
    \ddtau\H(\rho,w) + \D(\rho,w)  = -h_\partial m|_0^\ell
\end{align}
with dissipation functional $\D(\rho,w) := \int_0^\ell a \gamma\rho|w|^3\ dx  \geq 0$.
The system energy thus only changes by dissipation and energy flow across the boundary.
\end{lemma}
\begin{proof}
By formal differentiation, the identities \eqref{eq:energy-varder} and Lemma~\ref{lem:var}, we immediately get
\begin{align*}
    \ddtau\H(\rho,w) =&\ \la \delta_\rho\H,\ptau\rho\ra + \la \delta_w\H,\ptau w\ra
    = \la ah,\ptau\rho\ra + \la\eps^2m,\ptau w\ra\\ 
    =& -\la \dx m,h\ra + \la h,\dx m\ra - h_\partial m \vert_0^\ell - \la\gamma |w| w, m\ra
    = - h_\partial m \vert_0^\ell - \int_0^\ell a \gamma\rho|w|^3 dx,
\end{align*}
which already proves the required identity. \end{proof}

\begin{remark} \label{rem:pH}
Let us emphasize that the energy-identity \eqref{eq:energy} follows immediately from the variational identities \eqref{eq:var1}--\eqref{eq:var2} and the constitutive relations \eqref{eq:energy-varder} connecting the state and co-state variables.  
The particular form of the weak formulation also reveals the underlying \emph{port-Hamiltonian} structure, which can be preserved under Galerkin projection; see e.g. \cite{Egger19} for details. 
This will be the rationale behind our structure-preserving discretization strategy.
\end{remark}

For later reference, let us state a further important property of the energy functional. 
For ease of notation, we abbreviate $u=(\rho,w)$ and introduce the $\eps$-weighted norms
\begin{align}\label{eps-norm}
\|u\|_\eps^2 := \|\rho\|_{L^2(0,\ell)}^2 + \eps^2\|w\|_{L^2(0,\ell)}^2,
\qquad
    \|u\|_{\eps,\infty} := \|\rho\|_{L^\infty(0,\ell)} + \eps\|w\|_{L^\infty(0,\ell)},
\end{align}
which are well defined for all bounded measurable functions.
\begin{lemma} \label{lemma:H}
Let assumptions \ref{A1}--\ref{A2} hold. Then the energy functional $\H$ is well-defined, smooth, and uniformly convex on the set
\begin{align}\label{eq:bounded}
\AS:=\{(\rho,w) \in L^\infty(0,\ell)^2 : 2 \ubar \rho/3 \le \rho \le 3\bar \rho/2, \ -3\bar w/2 \le w \le 3\bar w/2\}
\end{align}
of \emph{admissible states} with respect to the weighted norm $\|\cdot\|_\eps$, i.e.,
\begin{align} \label{eq:estH}
    \H(u) - \H(\hat u) - \langle \H'(\hat u),u-\hat u\rangle \ge \tfrac{\alpha}{2} \|u-\hat u\|_\eps^2
\end{align}
for all $u,\hat u \in \AS$ and with some constant $\alpha>0$ independent of the parameter $\eps$.
\end{lemma}
\begin{remark}
The bounds in the definition of the set $\AS$ are slightly weaker than those in assumption \ref{A3}. Hence, sufficiently good approximations of a subsonic bounded state solution will therefore be admissible in the above sense, which will be used in the following.
\end{remark}
\begin{proof}
We show that the Hessian of the integrand in \eqref{eq:H}, which is given by
\begin{align} \label{eq:hessian}
    \H''(u) := \begin{pmatrix} \delta_{\rho\rho} \H & \delta_{\rho w} \H \\ \delta_{w \rho} \H & \delta_{ww}\H\end{pmatrix} 
    = \begin{pmatrix} a P''(\rho) & a \eps^2 w \\ a \eps^2 w & a \eps^2 \rho\end{pmatrix},
\end{align}
is positive definite for any $u=(\rho,w)$ satisfying the bounds in \eqref{eq:bounded}. To do so, we multiply $\H''(u)$ from left and right by $z=(x,y)$ and see that 
\begin{align*}
    z^\top \H''(u) z 
    &= aP''(\rho)x^2 + 2a\eps^2 w x y + a\eps^2\rho y^2\\
    &\ge  aP''(\rho)x^2 - \tfrac{4}{3} a \eps^2 \tfrac{9}{4}\tfrac{\bar w^2}{\rho} x^2 - \tfrac{3\rho}{4} a \eps^2 y^2 + a\eps^2\rho y^2\\
    &\ge a (1-\tfrac{3}{4}) P''(\rho) x^2 + a (1- \tfrac{3}{4}) \eps^2 \rho y^2 
    \ge \tfrac{\ubar a}{4} c_P x^2 + \tfrac{\ubar a}{4} \eps^2 \tfrac{2}{3}\ubar \rho y^2,
\end{align*}
where we used Young's inequality as well as \ref{A1}--\ref{A2}, which yield $P''(\rho) \ge c_P$ for some constant $c_P>0$.
The estimate \eqref{eq:estH} then follows with $\alpha=\min\{\frac{\ubar a c_P}{4},\frac{\ubar a\ubar \rho}{6}\}$ by Taylor expansion and integration over the spatial domain.
\end{proof}

\section{Structure-preserving discretization}
\label{sec:disc}

For the numerical approximation of the system \eqref{eq:gas1}--\eqref{eq:gas4}, we use a combination of a mixed finite-element method in space and the implicit Euler method in time. 
Let $x_i = i h$, $0 \le i \le M$ with $h=\ell/M$ be the grid points and $\T_h = \{T_i : 1 \le i \le M\}$ denote the corresponding mesh consisting of elements $T_i=[x_{i-1},x_i]$. 
We denote by 
\begin{align*}
    Q_h=\mathcal{P}_0(\T_h),\qquad R_h=\mathcal{P}_1(\T_h)\cap H^1(0,\ell)
\end{align*}
the spaces of piecewise constant and continuous piecewise linear functions over the mesh~$\T_h$. 
We further designate by $\Pi_h : L^2(0,\ell) \to P_0(\T_h)$ and $I_h : H^1(0,\ell) \to P_1(\T_h)$ the $L^2$-orthogonal projection and the piecewise linear interpolation operator, respectively.
We next define discrete time steps $\tau^n=n\, \Dtau$, $ n=0,\ldots,N$, with $\Dtau = \tau_{max}/N$ fixed for simplicity, and finally denote by $\dtau u^n = \tfrac{1}{\Dtau}(u^n-u^{n-1})$ the backwards difference quotient. 

\subsection{Definition of the discretization scheme}

For the numerical approximation of the system \eqref{eq:gas1}--\eqref{eq:gas4} on a single pipe, we consider the  following method.
\begin{problem}[Fully discrete scheme] \label{prob:fds}
Let $\rho_h^0 = \Pi_h \rho(0)$ and $m_h^0 = I_h m(0)$ be given.
Then for all $1 \le n \le N$, find $\rho_h^n\in Q_h,\ m_h^n\in R_h$ such that 
\begin{align}
	\la a \dtau \rho_h^n, q_h\ra + \la \dx m_h^n, q_h\ra =&\ 0\qquad \forall q_h\in Q_h,  \label{eq:discr1} \\
	\la \eps^2 \dtau w_h^n, r_h\ra - \la h_h^n, \dx r_h\ra + h_{\partial}^n r_h \vert_0^\ell
	+ \la\gamma |w_h^n|w_h^n, r_h\ra =&\ 0\qquad \forall r_h\in R_h, \label{eq:discr2}
\end{align}
with $h_\partial^n:=h_\partial(\tau_n)$, $w_h^n:=\tfrac{m_h^n}{a\rho_h^n}$, and $h_h^n:=\tfrac{\eps^2(m_h^n)^2}{2a^2(\rho_h^n)^2}+P'(\rho_h^n)$ introduced for abbreviation.
\end{problem}

\subsection{Basic properties of the discretization scheme}

The following result summarizes some of the basic properties of the numerical scheme introduced above.
\begin{lemma} \label{lem:varh}
Let $(\rho_h^{n-1},m_h^{n-1}) \in Q_h \times R_h$ be given with 
$2\ubar \rho/3 < \rho_h^{n-1} < 3 \bar \rho/2$ and $|w_h^{n-1}| < 3\bar w/2$. 
Then for any $0 < \Dtau \le \Dtau_0$ sufficiently small, the system \eqref{eq:discr1}--\eqref{eq:discr2} has a unique solution $(\rho_h^n,m_h^n) \in Q_h \times R_h$ in a small neighborhood of $(\rho_h^{n-1},m_h^{n-1})$ such that
$$
2\ubar \rho/3 < \rho_h^n < 3 \bar \rho/2 
\qquad \text{and} \qquad 
|w_h^n| < 3 \bar w/2.
$$
Moreover, any solution $(\rho_h^n,m_h^n)$ of \eqref{eq:discr1}--\eqref{eq:discr2} with the above bounds further satisfies 
\begin{align*}
    \dtau \H(\rho_h^n,w_h^n) + \D(\rho_h^n,w_h^n) \le -h_\partial^n m_h^n |_0^\ell.
\end{align*}
\end{lemma}
\begin{proof}
The first claim follows by a homotopy argument:
For $\Dtau=0$ the existence of a unique solution $\rho_h^n=\rho_h^{n-1}$ and $m_h^n=m_h^{n-1}$ is trivial. Since the Jacobian of the nonlinear system is regular for $\Dtau=0$ and depends continuously on the time-step, existence and local uniqueness of the solution follows by the  implicit function theorem.
To show the second claim, we note that $ f(u^n) - f(u^{n-1}) \le  f'(u^n) (u^n - u^{n-1}) $ for any smooth convex function $f$. 
From this and Lemma~\ref{lemma:H}, we conclude that
\begin{align*}
\dtau \H(\rho_h^n,w_h^n) 
&\le \la \delta_\rho \H(\rho_h^n,w_h^n), \dtau \rho_h^n \ra + \la \delta_w \H(\rho_h^n,w_h^n),\dtau w_h^n\ra \\
&= \la h_h^n, a \dtau \rho_h^n\ra + \la m_h^n, \eps^2  \dtau w_h^n\ra
= \la \Pi_h h_h^n, a \dtau \rho_h^n\ra + \la m_h^n, \eps^2  \dtau w_h^n\ra,
\end{align*}
where we used that $\H$ is convex on the set $\AS$ of admissible states, the relations between the discrete state and co-state variables announced in the lemma, as well as the orthogonality of the $L^2$-projection $\Pi_h$. 
%
The second assertion then follows by employing the identities  \eqref{eq:discr1}--\eqref{eq:discr2} with test functions $q_h=\Pi_h h_h^n$ and $r_h=m_h^n$, respectively.
\end{proof}

\begin{remark} \label{rem:AS}
By induction, one can see that for $\Dtau$ sufficiently small, the numerical solution will stay in the set $\AS$ of admissible states at least for a couple of time steps, if the initial values satisfy the bounds of the lemma.
We may therefore assume that 
\begin{enumerate}[label=(A\arabic*h)] \itemindent1.2em
    \setcounter{enumi}{2}
    \item \label{A3h}
    \qquad $2\ubar \rho/3 \le \rho_h^n \le 3\bar \rho/2$ \quad and \quad $|w_h^n| \le  3\bar w/2$ \qquad for all $0 \le n \le N^*$ 
\end{enumerate}
up to a certain index $N^*$ which may be smaller than $N=\tau_{max}/\Dtau$ in general. 
The index $N^*$ will however increase when $\Dtau$ and $h$ decrease. 
We will later see that $N^*=N$ if $\Dtau$ and $h$ are sufficiently small. 
The second assertion of the lemma shows that the proposed scheme inherits the energy-dissipation property of the continuous problem, i.e., it is a \emph{structure-preserving discretization scheme}.
\end{remark}

\subsection{Uniform convergence}

We are now in the position to state and prove our main result, which is concerned with the convergence of the discretization scheme above.
\begin{theorem}\label{thm:error_estimate}
Let \ref{A1}--\ref{A2} hold and $(\rho,w)$ denote a classical solution of \eqref{eq:gas1}--\eqref{eq:gas4} satisfying assumptions \ref{A3}--\ref{A4}. 
Then for any $0 < \Dtau \le \Dtau_0$ sufficiently small and $h\approx\Dtau$, Problem~\ref{prob:fds} has a unique discrete admissible state solution $(\rho_h^n,m_h^n)_{0 \le n \le N}$ which further satisfies 
\begin{align*}
	\|\rho(\tau^n)-\rho_h^n\|_{L^2(0,\ell)}^2 + \eps^2\|m(\tau^n)-m_h^n\|_{L^2(0,\ell)}^2 +  \sum_{k=1}^n \Dtau  \|m(\tau^k)-m_h^k\|_{L^3(0,\ell)}^3
	\leq C(\Dtau^2+h^2).
\end{align*}
The constants $C$ and $\Dtau_0$ in this assertion can be chosen independent of $\eps$.
\end{theorem}

The proof of this theorem basically relies on stability of the discrete problem and projection error estimates. Since we face a nonlinear problem, the \emph{relative energy} technique will be used for the former. 
A main subtlety here is to carefully track the dependence on the parameter $\eps$ in all estimates.
Details will be discussed in the next section. 

\begin{remark}
From the uniform error estimates, one can deduce that $(\rho_h^n,w_h^n) \in \AS$ for all $0 \le n \le N$ if the meshing parameters are sufficiently small, which is why we call the discrete solution \emph{admissible}.
Further note that the theorem yields uniform convergence estimates for all $\eps\ge 0$, in particular for $\eps=0$, which represents the parabolic limit problem. 
Hence, the proposed method and our estimates are \emph{asymptotic preserving}, i.e., for $\eps=0$ we obtain a viable discretization scheme for the parabolic limit problem with order optimal convergence rates. 
It will become clear later on, that the same convergence rates as for $m$ also hold for the velocity $w$, which plays an essential role in the stability analysis.
\end{remark}

\section{Proof of Theorem~\ref{thm:error_estimate}}
\label{sec:proof}

Before going into the details, let us briefly discuss the main arguments of our analysis. 
We define projections $\hat \rho_h(\tau) = \Pi_h \rho(\tau)$ and $\hat m_h(\tau) = I_h m(\tau)$ and abbreviate 
\begin{align*}
\hat \rho_h^n = \hat \rho_h(\tau^n) = \Pi_h \rho(\tau^n)
\qquad \text{and} \qquad 
\hat m_h^n = \hat m_h(\tau^n) = I_h m(\tau^n).
\end{align*}
By the triangle inequality, we can then decompose the error by
\begin{align}
    \|\rho(\tau^n)-\rho_h^n\|_{L^p(0,\ell)} &\le 
    \|\rho(\tau^n)-\hat\rho_h^n\|_{L^p(0,\ell)}
    + \|\hat\rho_h^n-\rho_h^n\|_{L^p(0,\ell)},\label{eq:split1}\\
    \|m(\tau^n)-m_h^n\|_{L^p(0,\ell)} &\le 
    \|m(\tau^n)-\hat m_h^n\|_{L^p(0,\ell)}
    + \|\hat m_h^n-m_h^n\|_{L^p(0,\ell)},\label{eq:split2}
\end{align}
into \emph{projection errors} and \emph{discrete error components}. 
The former can be estimated by standard arguments, and the main difficulty therefore is to show the respective bounds for the discrete error components.
By inserting the projections into the weak formulation of the problem, we can define residuals $\resq^n\in Q_h$, $\resr^n\in R_h$, given by
\begin{align}
	\la a \dtau \hat\rho_h^n, q_h\ra + \la \dx \hat m_h^n, q_h\ra &=: \la \resq^n,q_h\ra\qquad \forall q_h\in Q_h,  \label{eq:pert1}  \\
	\la \eps^2 \dtau \hat w_h^n, r_h\ra - \la \hat h_h^n, \dx r_h\ra + \hat h_{\partial}^n r_h \vert_0^\ell
	+ \la \gamma |\hat w_h^n|\hat w_h^n, r_h\ra &=: \la \resr^n,r_h\ra\qquad \forall r_h\in R_h. \label{eq:pert2} 
\end{align}
Like before, we use $\hat h_\partial^n:=h_\partial(\tau_n)$, $\hat w_h^n:=\tfrac{\hat m_h^n}{a\hat\rho_h^n}$, and $\hat h_h^n:=\tfrac{\eps^2(\hat m_h^n)^2}{2a^2(\hat\rho_h^n)^2}+P'(\hat\rho_h^n)$ for abbreviation. 
The projections can thus be understood as solutions of a perturbed discrete problem. 
To estimate the difference between the discrete solution and the projections, we will utilize \emph{relative energy estimates}, for which we require $(\rho_h^n,w_h^n)$ and $(\hat \rho_h^n,\hat w_h^n)$ to be admissible in the sense of condition~\ref{A3h}. Hence, our estimates will first only hold for the time steps $1 \le n \le N^*$, which may depend on the discretization parameters. 
As a consequence of the derived bounds, however, one can see that $N^*=N$ as soon as $\Dtau$ and $h$ are small enough. 

The remainder of this section is organized as follows: In Section~\ref{ssec:proj}, we state estimates for the projection errors
and show that $(\hat \rho_h^n,\hat w_h^n)$ is admissible for all $1 \le n \le N$ if $h$ is sufficiently small.  
In Section~\ref{ssec:rel}, we introduce and derive some properties of the relative energy. 
Section~\ref{ssec:reldiff} then contains a technical results concerning time differences of the relative energy and 
Section~\ref{ssec:aux} presents the relative energy estimates for the discrete error components. The proof of Theorem~\ref{thm:error_estimate} is finally completed in Section~\ref{ssec:proof}.

\subsection{Projection errors} \label{ssec:proj}

Let us start with summarizing some elementary properties of the projection operators $I_h : C([0,\ell]) \to P_1(\T_h) \cap H^1(0,\ell)$ and $\Pi_h : L^2(0,\ell) \to P_0(\T_h)$.
Also recall that $\T_h$ is a uniform mesh with elements $[x_{i-1},x_i]$ of size $h$ by assumption. 

\begin{lemma}\label{lemma:proj}
For any $z \in W^{1,p}(0,\ell)$, $1 \le p \le \infty$, there holds 
$$
\dx (I_h z) = \Pi_h (\dx z).
$$ 
Moreover, 
$\|\Pi_h z\|_{L^\infty(0,\ell)} \le \| z\|_{L^\infty(0,\ell)}$ and $\|I_h z\|_{L^\infty(0,\ell)} \le \|z\|_{L^\infty(0,\ell)}$, as well as
\begin{align*}
\|z - \Pi_h z\|_{L^p(0,\ell)} \le c h \|\dx z\|_{L^p(0,\ell)} 
\qquad \text{and} \qquad 
\|z - I_h z\|_{L^p(0,\ell)} \le c h \|\dx z\|_{L^p(0,\ell)}
\end{align*}
hold with generic constant $c$ independent of $z$, $p$, and $h$. By the local definition of the operators, all assertions also hold locally on every element $[x_{i-1},x_i]$.
\end{lemma}
\begin{proof}
The first assertion is known as \emph{commuting diagram} property and follows immediately from the fundamental theorem of calculus. 
The boundedness of the projection $\Pi_h$ follows by noting that $\Pi_h z(\xi_i)=z(\xi_i)$ for some $\xi_i \in [x_{i-1},x_i]$, and noting that $\Pi_h z$ is constant on every element $[x_{i-1},x_i]$. 
That for the interpolation follows with a similar argument.
The error estimates can for instance be found in \cite[Ch.~4]{BrennerScott_2008}. 
\end{proof}

As a direct consequence of the above estimates, we obtain the following assertions. 
\begin{lemma}\label{lemma:projerror}
Let \ref{A1}--\ref{A4} hold and $\hat \rho_h^n=\Pi_h \rho(\tau^n)$ and $\hat m_h^n = I_h m(\tau^n)$. Then
\begin{align*}
\|\rho(\tau^n) - \hat \rho_h^n\|^2_{L^2(0,\ell)} + \eps^2 \|m(\tau^n) - \hat m_h^n\|_{L^2(0,\ell)}^2 + \sum_{k=1}^n \Dtau \|m(\tau^k) - \hat m_h^k\|_{L^3(0,\ell)}^3 \le C h^2     
\end{align*}
with constant $C$ only depending on the bounds in the assumptions. 
Moreover, 
$\ubar \rho \le \hat \rho_h^n \le \bar \rho$
and for any $0 < h \le h_0$ sufficiently small, we have $-3\bar w/2 \le \hat w_h^n \le 3\bar w/2$
where $\hat w_h^n = \frac{\hat m_h^n}{a \hat \rho_h^n}$. 
\end{lemma}
\begin{proof}
The error estimate is a direct consequence of the bounds stated in the previous lemma and the regularity of the solution provided by condition \ref{A4}. 
From the properties of the projection operators $\Pi_h$ and $I_h$ and assumption~\ref{A3}, one can see that $\ubar \rho \le \hat \rho_h^n \le \bar \rho$, which implies the pointwise bounds for the density, as well as $|m|\le \bar \rho \bar w$ and $|\hat m_h^n| \le \bar \rho \bar w$.
For the velocity $w=\frac{m}{a\rho}$ and its discrete counter part $\hat w_h=\frac{\hat m_h}{a \hat \rho_h}$, we may expand
\begin{align*}
w - \hat w_h = \tfrac{1}{a\rho} (m - \hat m_h) + \tfrac{\hat m_h}{a \rho \hat \rho_h} (\hat \rho_h - \rho).
\end{align*}
By the uniform bounds for density and mass flux and the approximation error estimates of the previous lemma, we hence conclude that
\begin{align*}
    \|w(\tau^n) - \hat w_h^n\|_{L^\infty(0,\ell)} 
    \le \tfrac{1}{a\ubar \rho} \|m(\tau^n) - \hat m_h^n \|_{L^\infty(0,\ell)} + \tfrac{\bar\rho\bar w}{\ubar\rho^2} \|\rho(\tau^n) - \hat \rho_h^n\|_{L^\infty(0,\ell)} 
    \le C h,
\end{align*}
with some constant $C$ depending only on the bounds in the assumptions. For $h$ small enough, we have $C h \le \bar w/3$, which yields the pointwise bounds for $\hat w_h^n$.
\end{proof}
The projections of smooth bounded subsonic state solutions therefore always lie in the set $\AS$ of admissible states, as soon as the mesh size $h$ is sufficiently small.
It is further clear from the proof, that one could simply replace $m$ by $w$ in the error estimate 
to obtain corresponding approximation error bounds for the velocity.

\subsection{Relative energy} \label{ssec:rel}

In order to measure the distance between two (approximate) solutions, we will utilize the concept of \emph{relative energy}; see \cite{Dafermos05}. 
For ease of notation, we abbreviate $u=(\rho,w)$, $\hat u=(\hat\rho,\hat w)$ in the following, and recall the definition
\begin{align}\label{def:relen}
	\H(u\mid\hat{u}) 
	:= \H(u) - \H(\hat{u})
	- \langle\H'(\hat{u}), u-\hat{u}\rangle
\end{align}
of the relative energy.
Since $\H$ is strictly convex on the set of admissible states, see Lemma~\ref{lemma:H}, the relative energy $\H(u\mid \hat{u})$ is positive on the set $\AS$ and defines a distance measure there, which is equivalent to the $\eps$-weighted $L^2$-norm defined in \eqref{eps-norm}.
\begin{lemma}\label{lemma:relen}
Let assumptions \ref{A1}--\ref{A2} hold and 
$u,\hat u \in \AS$ be admissible. Then 
\begin{align}\label{eq:relen:1}
c_0 \|u-\hat u\|_\eps^2
\le \H(u\mid \hat u) \le C_0 \|u-\hat u\|_\eps^2,
\end{align}
and for all $x \in L^\infty(0,\ell)^2 ,\ y\in L^2(0,\ell)^2$, one has 
\begin{align}
    \langle (\H''(u)-\H''(\hat u))\, x,y\rangle &\le C\|u-\hat u\|_{\eps}\|x\|_{\eps,\infty}\|y\|_{\eps}. \label{ddH:cond:2}
\end{align}
The constants $c_0,\, C_0,\, C$ in these estimates only depend on the bounds in \ref{A1}--\ref{A2}.
\end{lemma}
\begin{proof}
Define $F(s) := \H(s u + (1-s) \hat u)$. Then by Taylor's theorem and the chain rule
\begin{align*}
\H(u\mid \hat u) &= 
\H(u) - \H(\hat u)
- \langle\H'(\hat u), u-\hat u \rangle 
= F(1) - F(0) - F'(0) \\
&= \tfrac{1}{2} F''(s^\ast) 
=\tfrac{1}{2}\langle \H''(u^\ast)(u-\hat u),u-\hat u\rangle
\end{align*}
for some $0 < s^\ast < 1$ and $u^\ast:=s^\ast u + (1-s^\ast) \hat u$ satisfying the bounds of the lemma, since the set $\AS$ is convex.
The lower bound then is a direct consequence of $\H$ being strictly convex for bounded subsonic states w.r.t. $\|\cdot\|_\eps$, see Lemma~\ref{lemma:H}. 
From the formula for the Hessian in \eqref{eq:hessian}, one can further see that
\begin{align*}
\langle\H''(u^\ast) (u-\hat{u}),u-\hat{u}\rangle
\le C_0 \|u-\hat u\|_\eps^2,
\end{align*}
where we used that $P$ is smooth, with some upper bound $C_P$ for $P''(\rho)$, as well as Young's inequality and the bounds in the definition of the set $\AS$; this proves the first assertion. 
For functions $u,\hat u\in \AS$ and $x\in L^\infty(0,\ell)^2,\,y\in L^2(0,\ell)^2$, one can further see that
\begin{align*}
    \langle (\H''(u)-&H''(\hat u))x,y\rangle = (a(P''(\rho) - P''(\hat\rho))x_1,y_1) + (a\eps^2(w-\hat w)x_1,y_2)\\
    &\qquad\qquad\qquad\qquad\qquad\ + (a\eps^2(w-\hat w)x_2,y_1) + (a\eps^2 (\rho-\hat\rho)x_2,y_2)\\
    \le&\ C a \|\rho-\hat\rho\|_{L^2}\|x_1\|_{L^\infty}\|y_1\|_{L^2}
    + a\|\eps w-\eps \hat w\|_{L^2}\|x_1\|_{L^\infty}\|\eps y_2\|_{L^2}\\
    &\qquad\qquad \ + a\|\eps w-\eps \hat w\|_{L^2}\|\eps x_2\|_{L^\infty}\|y_1\|_{L^2}
    + a \|\rho-\hat\rho\|_{L^2}\|\eps x_2\|_{L^\infty}\|\eps y_2\|_{L^2},
\end{align*}
where we used that the third derivative of $P$ is bounded on the set $\AS$. The terms in the estimate can further be bounded by $C' \|u-\hat u\|_{\eps}\|x\|_{\eps,\infty}\|y\|_{\eps}$ by definition of the norms.
\end{proof}

\subsection{Time differences of the relative energy} \label{ssec:reldiff}

In order to measure the distance between two discrete functions recursively, we will utilize the following technical result. 

\begin{lemma}\label{lemma:ddH} 
Let \ref{A1}--\ref{A2} hold and assume that $u^k$, $\hat u^k\in \AS$ for $k \in \{n-1,n\}$. Then 
\begin{align} \label{eq:ddH0}
    \dtau \H(u^n|\hat u^n) \leq&\ \langle \H'(u^n) - \H'(\hat u^n) - \H''(\hat u^n)(u^n-\hat u^n),\dtau\hat u^n \rangle\\
    & + \langle \H'(u^n)- \H'(\hat u^n),\dtau u^n - \dtau\hat u^n\rangle\nonumber\\
    &+ C \|\dtau \hat u^n\|_{\eps,\infty} \big( \H(u^n|\hat u^n) + \H(u^{n-1}|\hat u^{n-1}) \big) + C'\|\dtau \hat u^n\|_{\eps,\infty}\|\hat u^n - \hat u^{n-1}\|_{\eps}^2\nonumber
\end{align}
with constants $C,\ C'$ independent of $\eps$ and $\Dtau$ as well as the functions $u^k$ and $\hat u^k$. 
\end{lemma}

\begin{proof}
Using Taylor's theorem and rearranging the terms, we observe that 
\begin{align*}
	\dtau \H(u^n|\hat u^n)
	=\tfrac{1}{\Dtau}\big(\H(&u^n) - \H(u^{n-1}) - \H(\hat u^n) + \H(\hat u^{n-1})\\
	&- \langle \H'(\hat u^n),u^n-\hat u^n\rangle + \langle \H'(\hat u^{n-1}),u^{n-1}-\hat u^{n-1} \rangle  \big)\\
	=\langle \H'(u^n), \dtau u^n\rangle -\tfrac{\Dtau}{2} 
	&\langle \H''(u^\ast) \, \dtau u^n, \dtau u^n \rangle
	-\langle \H'(\hat u^n), \dtau \hat u^n\rangle +\tfrac{\Dtau}{2} 
	\langle \H''(\hat u^\ast) \, \dtau \hat u^n, \dtau \hat u^n \rangle\\
	&-\langle \H'(\hat u^n), \dtau  u^n -\dtau \hat u^n\rangle
	- \tfrac{1}{\Dtau}\langle  \H'(\hat u^n)- \H'(\hat u^{n-1}), u^{n-1} - \hat u^{n-1} \rangle
\end{align*}
with intermediate values 
$u^\ast,\hat u^\ast$ lying on the lines between $u^n$ and $ u^{n-1}$ and $\hatn u$ and $\hat u_h^{n-1}$, respectively; 
in particular, $u^\ast,\hat u^\ast \in \AS$.
By suitably adding and subtracting terms of the form $\langle \H'(u^n),\dtau\hat u^n\rangle$, $ \langle\dtau\H'(\hat u^n),u^n-\hat u^n\rangle$ and $\langle\H''(\hat u^n)(u^n-\hat u^n), \dtau \hat u^n\rangle$, one can further see that
\begin{align}\label{eq:derivH:1}
	\dtau \H(u^n\mid \hat u^n) 
	=&\ \langle \H'(u^n)-\H'(\hat u^n), \dtau u^n-\dtau \hat u^n\rangle\\
	&\ + \langle\H'(u^n)-\H'(\hat u^n) - \H''(\hat u^n)(u^n-\hat u^n), \dtau \hat u^n\rangle\nonumber\\
	&\ \ - \tfrac{\Dtau}{2} \langle\H''(u^\ast) \dtau u^n, \dtau u^n\rangle
	+ \tfrac{\Dtau}{2} \langle\H''(\hat u^\ast)\,  \dtau \hat u^n, \dtau \hat u^n\rangle\nonumber\\
	&\ \ \ -\langle \dtau \H'(\hat u^n) - \H''(\hat u^n)\, \dtau \hat u^n, u^n-\hat u^n \rangle
	+ \Dtau \langle \dtau\H'(\hat u^n),\dtau u^n - \dtau\hat u^n\rangle.\nonumber
\end{align}
The terms in the first two lines already appear in the final estimate, and by Taylor's theorem, the last two lines can be transformed into
\begin{align*}
(*)=- \tfrac{\Dtau}{2}& \langle\H''(u^\ast) \dtau u^n, \dtau u^n\rangle
	+ \tfrac{\Dtau}{2} \langle\H''(\hat u^\ast)\,  \dtau \hat u^n, \dtau \hat u^n\rangle\nonumber\\
	&-\langle (\H''(\hat u^{\ast\ast})- \H''(\hat u^n))\, \dtau \hat u^n, u^n-\hat u^n \rangle
	+\Dtau\langle \H''(\hat u^{\ast\ast\ast})\dtau\hat u^n,\dtau u^n - \dtau\hat u^n\rangle
\end{align*}
with intermediate values  
$\hat u^{\ast\ast},\ \hat u^{\ast\ast\ast} \in \AS$.
After rearranging the terms and expanding by 
$\Dtau\langle\H''(u^\ast)\dtau\hat u^n,\dtau u^n\rangle$ and $  \tfrac{\Dtau}{2}\langle \H''(u^\ast)\dtau\hat u^n,\dtau\hat u^n\rangle$, we obtain
\begin{align*}
	(*)=- \tfrac{\Dtau}{2}& \langle\H''(u^\ast) \dtau u^n, \dtau u^n\rangle
	+\Dtau  \langle \H''(u^\ast) \dtau \hat u^n, \dtau u^n\rangle 
	- \tfrac{\Dtau}{2} \langle \H''(u^\ast) \dtau\hat u^n, \dtau \hat u^n\rangle\\
	&+\Dtau \langle \big( \H''(\hat u^{\ast\ast\ast})-\H''(u^\ast)\big)\dtau \hat u^n, \dtau u^n\rangle + \tfrac{\Dtau}{2}\langle \big(\H''(\hat u^\ast) - \H''(\hat u^{\ast\ast\ast})\big) \dtau\hat u^n , \dtau\hat u^n \rangle\\
	&\quad +\tfrac{\Dtau}{2} \langle \big(\H'' (u^\ast) - \H''(\hat u^{\ast\ast\ast})\big) \dtau \hat u^n , \dtau \hat u^n\rangle
    -\langle (\H''(\hat u^{\ast\ast})- \H''(\hat u^n))\, \dtau \hat u^n, u^n-\hat u^n \rangle.
\end{align*}
Since $\H$ is convex on the set $\AS$, the first line is non-positive, i.e.,
\begin{align*}
	-\tfrac{\Dtau}{2} \langle \H''(u^\ast) (\dtau u^n-\dtau\hat u^n), \dtau u^n-\dtau\hat u^n \rangle \le 0.
\end{align*}
Inequality \eqref{ddH:cond:2} then further allows us to estimate the remaining four terms such that
\begin{align*}
    (*) \le C\big(\|\dtau& \hat u^n\|_{\eps,\infty}\|u^\ast - \hat u^{\ast\ast}\|_{\eps} \|u^n - u^{n-1}\|_{\eps}
    + \|\dtau \hat u^n\|_{\eps,\infty}\|\hat u^{\ast} - \hat u^{\ast\ast\ast}\|_{\eps} \|\hat u^n - \hat u^{n-1}\|_{\eps}\\
    &+ \|\dtau \hat u^n\|_{\eps,\infty}\|u^{\ast} - \hat u^{\ast\ast\ast}\|_{\eps} \|\hat u^n - \hat u^{n-1}\|_{\eps}
    + \|\dtau \hat u^n\|_{\eps,\infty}\|\hat u^n - \hat u^{\ast\ast} \|_{\eps} \|u^n - \hat u^n\|_{\eps}\big).
\end{align*}
By elementary manipulations, one can see that
\begin{align*}
    \|u^n-u^{n-1}\|_{\eps} &\le \|u^n-\hat u^n\|_{\eps} + \|\hat u^n-\hat u^{n-1}\|_{\eps}  + \|u^{n-1}-\hat u^{n-1}\|_{\eps},\\
    \|u^\ast - \hat u^{\ast\ast}\|_{\eps} 
    &\le \|u^n-\hat u^n\|_{\eps} + \|\hat u^n-\hat u^{n-1}\|_{\eps}  + \|u^{n-1}-\hat u^{n-1}\|_{\eps},\\
    \|\hat u^{\ast} - \hat u^{\ast\ast}\|_{\eps} &\le \|\hat u^n - \hat u^{n-1}\|_{\eps},\qquad \|\hat u^n - \hat u^{\ast\ast\ast}\|_{\eps} \le \|\hat u^n - \hat u^{n-1}\|_{\eps},
    \end{align*}
which together with Lemma \ref{lemma:relen} and the previous calculations yields the assertion. \end{proof}

\begin{remark} 
On the continuous level, the time derivative of the relative energy is given by
\begin{align*}
    \frac{d}{d\tau} \H(u|\hat u) &= \langle \H'(u) - \H'(\hat u) - \H''(\hat u)(u-\hat u),\ptau\hat u \rangle
     + \langle \H'(u)- \H'(\hat u),\ptau u - \ptau\hat u\rangle.
\end{align*}
The terms in the last line of the estimate in Lemma \ref{lemma:ddH} hence are perturbations that are caused by the time discretization.
Further note that the result is only based on smoothness and convexity of the energy functional and thus independent of the particular problem.
\end{remark}

\subsection{Relative energy estimates}\label{ssec:aux}

We now turn our attention to the discrete error. 
By carefully estimating the terms in the right  hand side of \eqref{eq:ddH0}, 
we show the following result.
\begin{lemma} \label{lem:ddH}
Let \ref{A1}--\ref{A4} hold and $0 < h \le h_0$, $0 < \Dtau \le \Dtau_0$ be sufficiently small, such that assumption \ref{A3h} is valid for all $n \le N^*$. Then 
\begin{align}\label{eq:ddH}
    \dtau\H(u_h^n|\hatn u) \le&\ C\H(u_h^n|\hatn u) + C'\H(u_h^{n-1}|\hat u_h^{n-1}) + \tfrac{1}{\Dtau}(h-\hat h_h)(\rho_h-\hat\rho_h)|_{\tau_{n-1}}^{\tau_n}\\
    &\qquad+C'' (\Dtau^2+h^2) 
    -\tfrac{1}{2}\D(u_h^n|\hatn u)\nonumber
\end{align}
holds for all $n \le N^*$ with relative dissipation functional defined by
 \begin{align}
    \D (\hn u\mid \hatn u) := \int_0^\ell \tfrac{1}{4} \gamma a \hatn\rho  |\hn w -\hatn w|^2 (|\hatn w |+ |\hn w|)\, dx
    \ge \tfrac{1}{4} a\ubar \gamma \ubar\rho \|w_h^n-\hatn w\|_{L^3(0,\ell)}^3\ge 0.	\label{def:D}
\end{align} 
The constants $h_0$, $\Dtau_0$, $C$, $C'$ and $C''$ depend only on the bounds in the assumptions. 
\end{lemma}
A discrete Gronwall estimate and Lemma \ref{lemma:relen} then already yield the following estimate for the discrete error components for time steps $n \le N^*$, i.e,
\begin{align} \label{eq:rhowdisc}
	\|\hat \rho_h^n-\rho_h^n\|_{L^2(0,\ell)}^2 + \eps^2\|\hat w_h^n-w_h^n\|_{L^2(0,\ell)}^2 +   \sum_{k=1}^n \Dtau \|\hat w_h^k-w_h^k\|_{L^3(0,\ell)}^3
	\leq C(\Dtau^2+h^2).    
\end{align}
Using that $\hat m_h^n = a \hat \rho_h^n \hat w_h^n$ and $m_h^n =a \rho_h^n w_h^n$ by definition, one can show the same bounds with $w$ replaced by $m$, which already leads to the estimate of Theorem~\ref{thm:error_estimate}. 
\subsection*{Proof of Lemma~\ref{lem:ddH}}
Inequality~\eqref{eq:ddH} follows directly from \eqref{eq:ddH0} by appropriately estimating the terms in the three lines of the right hand side of the latter. 
For the remainder of this section, we assume \ref{A1}--\ref{A4} to be true.

\subsection*{Step~1.}

By the properties of the projections and assumption \ref{A4}, we can bound 
\begin{align*}
\|\dtau \hatn u\|_{\eps,\infty} 
\le \|\ptau u\|_{L^\infty(\tau_{n-1},\tau_n;\,\|\cdot\|_{\eps,\infty})} \le \bar C,
\end{align*}
and by Taylor estimates and similar arguments, we further see that
\begin{align*}
\|\hatn u - \hat u_h^{n-1}\|_{\eps}^2 \le C\Dtau^2\|\ptau \hat u_h\|_{L^\infty(\tau_{n-1},\tau_n;\,\| \cdot \|_\eps)}^2 \le C'\Dtau^2\|\ptau u\|_{L^\infty(\tau_{n-1},\tau_n;\,\| \cdot \|_\eps)}^2 \le C''\Dtau^2.
\end{align*}
This already allows to estimate the two terms in the third line of \eqref{eq:ddH0} accordingly.

\subsection*{Step 2.} 
From the formulas for the derivatives of the energy functional, we see that 
\begin{align*}
   \H'(u_h^n) - \H'(\hatn u) - \H''(\hatn u)(u_h^n-\hatn u) 
   =\begin{pmatrix} aP'(\rho_h^n|\hatn\rho) + \tfrac{a\eps^2}{2}(w_h^n-\hatn w)^2\\ 
   a\eps^2(\rho_h^n-\hatn\rho)(w_h^n-\hatn w)\end{pmatrix}.
\end{align*}
By Taylor expansion and assumption \ref{A1}, we can further estimate
\begin{align*}
    | P'(\rho_h^n|\hatn\rho)| = | P'''(\rho_h^\ast)|(\rho_h^n-\hatn\rho)^2\le C \, |\rho_h^n-\hatn\rho|^2
\end{align*}
with intermediate value 
$\rho_h^\ast \in [\ubar \rho,\bar \rho]$. 
Using assumptions \ref{A1}--\ref{A3} as well as Hölder and Young inequalities, we further deduce that
\begin{align*}
    \langle \H'(u^n_h) - \H'(\hat u^n_h) - \H''(\hat u^n_h)(u^n_h-\hat u^n_h),\dtau\hat u^n_h \rangle \leq C'\|\dtau\hatn u\|_{\eps,\infty}\|u_h^n-\hatn u\|_{\eps}^2\le C''\H(u_h^n|\hatn u).
\end{align*}
Here we employed that $\|\dtau\hatn u\|_{\eps,\infty}$ is bounded, see Step~1, as well as Lemma \ref{lemma:relen}. 
This already yields the bound for the second term in the right hand side of \eqref{eq:ddH0}.

\subsection*{Step~3}

Bounding the 
second term on the right hand side of \eqref{eq:ddH0} turns out to be the most difficult task. 
By definition of the co-state variables, we see that
\begin{align*}
    \langle \dtau \hn u &- \dtau\hatn u, \H'(u_h^n)- \H'(\hatn u)\rangle\\
    &= (\dtau\rho_h^n - \dtau\hatn\rho,ah_h^n - a\hatn h) + (\dtau w_h^n - \dtau\hatn w,\eps^2m_h^n - \eps^2\hatn m) = (*).
\end{align*}
Since $(\rho_h^n,m_h^n)$ solves \eqref{eq:discr1}--\eqref{eq:discr2} and $(\hatn\rho,\hatn m)$ can be understood as solution of the perturbed system \eqref{eq:pert1}--\eqref{eq:pert2}, the above expression equals
\begin{align}\label{eq:step3:i}
    (*)&= -\, \la \dx m_h^n-\dx\hatn m,h_h^n-\hatn h\ra - \la \resq^n,h_h^n-\hatn h\ra
    +\la h_h^n-\hatn h,\dx m_h^n - \dx\hatn m\ra\nonumber\\
    &\qquad - \la\gamma|w_h^n|w_h^n - \gamma|\hatn w|\hatn w,m_h^n-\hatn m\ra - (h_{\partial}^n - h_{\partial}^n)(m_h^n-\hatn m)|_0^\ell - \la \resr^n,m_h^n-\hatn m\ra\nonumber\\
    &= -\, \big\la \gamma |\hn w|w_h^n - \gamma |\hatn w|\hatn w,\hn m - \hatn m\big\ra-\big\la \resq^n, \hn h - \hatn h\big\ra - \big\la \resr^n,\hn m - \hatn m\big\ra.
\end{align}
 By definition of $\resq^n$ in \eqref{eq:pert1} we have $\resq^n = a\dtau\hatn\rho + \dx\hatn m$ which can be tested with any $L^2$-function.
The following three lemmas provide the required estimates for the three terms in \eqref{eq:step3:i}. 
In the remainder of this section
we require assumptions \ref{A1}--\ref{A4} as well as condition \ref{A3h} to hold and $n \le N^*$. 

\begin{lemma} \label{lem:step3i}
The first term in \eqref{eq:step3:i} can be estimated by
\begin{align*}
    -\big\la \gamma |\hn w|w_h^n - \gamma |\hatn w|\hatn w,\hn m - \hatn m\big\ra
    \le  - \D(\hn u\mid \hatn u) + C \H(\hn u\mid \hatn u).
\end{align*}
\end{lemma}
\begin{proof}
We refer to \cite[Lemma~10]{EggerGiesselmann} for a detailed proof of this technical result.
\end{proof}

\begin{lemma}\label{lemma:res1}
The second term in \eqref{eq:step3:i} can be estimated by
\begin{align*}  
    -\la \resq^n,h_h^n-\hatn h\ra \le&\ C\Dtau^2 + C' \H(\hn u\mid\hatn u)
\end{align*}
with constants $C,\,C'$ that only depend on the bounds in the assumptions.
\end{lemma}

\begin{proof}
Due to the definition of $\resq^n$ in \eqref{eq:pert1} it holds that $\resq^n = a\dtau\hatn\rho + \dx\hatn m$. This identity can be tested with any $L^2$-function. We deduce
\begin{align*}
	-\big\la\resq^n, \hn h - \hatn h\big\ra &= - \la a \dtau \hatn \rho , \hn h - \hatn h\ra
	- \la \dx \hatn m, \hn h - \hatn h\ra\\
	&= -\la a \dtau \hatn \rho - a \ptau \hatn \rho , \hn h - \hatn h\ra
	- \la a \ptau \hatn \rho , \hn h - \hatn h \ra
	+ \la \dx \hatn m, \hn h - \hatn h\ra.
\end{align*}
The last two terms vanish due to the fact that $\hatn\rho,\,\dx\hatn m$ are piecewise constant in space,
which implies together with the definition of the projection $\Pi_h$ and $I_h$ and \eqref{eq:var1} that
\begin{align*}
	 \la a \ptau \hatn \rho,q\ra + \la \dx \hatn m,q \ra = 0\quad\forall q\in L^2(0,\ell).
\end{align*}
The first term can be estimated by Young's inequality, i.e., 
\begin{align*}
	-(a \dtau \hatn \rho - a \ptau \hatn \rho , \hn h - \hatn h) \le \tfrac{a^2}{2} \LL{\dtau \hatn \rho - \ptau \hatn \rho}
	+ \tfrac{1}{2} \LL{\hn h - \hatn h}.
\end{align*}
For the first term we deduce
\begin{align*}
	\LL{\dtau \hatn \rho - \ptau \hatn \rho}
	&\le \tfrac{\Dtau^2 }{4} \|\pptau \hat\rho_h\|_{L^\infty(\tau^{n-1}, \tau^n;L^2(0,\ell))}^2 \\
	&\le \tfrac{\Dtau^2 }{4} \|\pptau \rho\|_{L^\infty(\tau^{n-1}, \tau^n;L^2(0,\ell))}^2 \le C\Dtau^2,
\end{align*}
where the second estimate holds true due to the construction of the projection and the third due to~\ref{A4}.
In order to estimate the second term we use the bounds in \ref{A2}--\ref{A3} as well as the fact that the pressure potential $P$ is smooth by \ref{A1}. It then holds that
\begin{align*}
	|\hn h - \hatn h|
	= |\tfrac{\eps^2}{2} (|\hn w|^2-|\hatn w|^2) + P'(\hn\rho)-P'(\hatn\rho)|
	\le \bar\eps\bar w |\eps\hn w-\eps\hatn w| +C''|\hn\rho-\hatn\rho|.
\end{align*}
By Lemma \ref{lemma:relen} we conclude that $\LL{\hn h - \hatn h}$ can be estimated by $C'\H(\hn u\mid \hatn u)$.
\end{proof}

\begin{lemma}\label{lemma:II}
The third term in \eqref{eq:step3:i} can be estimated by
\begin{align*}  
    - \la \resr^n,m_h^n-\hatn m\ra \le&\  C (\Dtau^2+h^2) + C' \H(\hn u\mid\hatn u) + 
    C'' \H(u_h^{n-1}\mid \hat u_h^{n-1})\\
	&\qquad\qquad\qquad +\tfrac{1}{\Dtau } \la h-\hat h_h,\rho_h-\hat\rho_h\ra|_{\tau_{n-1}}^{\tau_n}
	+ \tfrac{1}{2}\D(\hn u\mid \hatn u)
\end{align*}
with 
constants $C$, $C'$ and $C''$ only depending on the bounds in the assumptions.
\end{lemma}

\begin{proof}
By consistency of the discrete problem, the exact solution $(\rho,m)$ of \eqref{eq:var1}--\eqref{eq:var2} satisfies 
\begin{align*}
    \la \eps^2\ptau w^n,r\ra - \la h^n,\dx r\ra + h_\partial^nr|_0^\ell + \la \gamma|w^n|w^n,r\ra  =0
\end{align*}
for all $r\in H^1(0,\ell)$, where we abbreviate $w^n = w(\tau_n),\ m^n = m(\tau_n),\ h^n = h(\tau_n)$. Adding this to $-\la \resr^n,m_h^n-\hatn m\ra$ with test function $r=m_h^n-\hatn{m}$ leads to 
\begin{align} \label{eq:res2est}
-\la \resr^n,m_h^n-\hatn m\ra = \la\eps^2\ptau w^n &- \eps^2\dtau\hatn w,m_h^n-\hatn m\ra - \la h^n - \hatn h,\dx m_h^n - \dx\hatn m\ra \\ 
&  + \la \gamma|w^n|w^n - \gamma|\hatn w|\hatn w,m_h^n-\hatn m\ra \nonumber
= (i) + (ii) + (iii).
\end{align}
Here we abbreviate $u^n=u(t^n)$ for continuous functions of time. 
In the following we consider the terms $(i)$--$(iii)$ separately, and use the projection error estimates from Lemma~\ref{lemma:proj} in order to get convergence rates.
Using Young's inequality we can estimate
\begin{align*}
    (i) = \la \eps^2\ptau w^n - \eps^2\dtau\hatn w,m_h^n-\hatn m\ra
    \le \tfrac{\eps^2}{2}\|\ptau w^n - \dtau\hatn w\|_{L^2}^2 + \tfrac{\eps^2}{2}\|m_h^n-\hatn m\|_{L^2}^2.
\end{align*}
By Lemma \ref{lemma:relen} the second term can be estimated by the relative energy, more precisely
\begin{align*}
  \tfrac{\eps^2}{2}\|m_h^n-\hatn m\|_{L^2}^2 \le&\ \tfrac{a^2}{2} \big(\eps^2\|\hatn\rho\|_{L^\infty}^2\|w_h^n-\hatn w\|_{L^2}^2 + \eps^2\|w_h^n\|_{L^\infty}^2\|\rho_h^n-\hatn\rho\|_{L^2}^2\big)
  \le C \H(\hn u\mid \hatn u),
\end{align*}
with constant $C$ only depending on the bounds in the assumptions.
The first term can be estimated by the triangle inequality, i.e.,
\begin{align*}
    \tfrac{\eps^2}{2}\|\ptau w^n - \dtau\hatn w\|_{L^2}^2
    \le \eps^2\|\ptau w^n - \ptau\hatn w\|_{L^2}^2 + \eps^2\|\ptau \hatn w - \dtau\hatn w\|_{L^2}^2.
\end{align*}
We now consider both terms separately: For the first one we deduce
\begin{align*}
    \eps\|\ptau w^n - &\ptau \hatn w\|_{L^2} 
    \le \eps\|\ptau \tfrac{m^n}{a\hat \rho_h^n} - \ptau \tfrac{\hat m_h^n}{a \hat \rho_h^n}\|_{L^2} + \eps \|\ptau \tfrac{m^n}{a \rho^n} - \ptau \tfrac{m^n}{a \hat \rho_h^n}\|_{L^2}\\
    &\le \tfrac{\eps}{a\ubar\rho}\|\ptau m^n-\ptau\hatn m\|_{L^2} + \tfrac{\eps}{a\ubar\rho^2}\|\ptau  m^n\|_{L^\infty} \|\rho^n-\hatn\rho\|_{L^2}+ \tfrac{\eps\bar \rho\bar w}{a\ubar\rho^2}\|\ptau\rho^n-\ptau\hatn\rho\|_{L^2}\\
    & \qquad \qquad \qquad \qquad 
    + \|\ptau\rho^n\|_{L^\infty}(\tfrac{\eps }{a\ubar\rho^2}\|m^n-\hatn m\|_{L^2} + \tfrac{3\eps \bar\rho^2\bar w}{a\ubar\rho^4}\|\rho^n-\hatn\rho\|_{L^2})
    \le C h
\end{align*}
where we used Lemma~\ref{lemma:proj}, the bounds of Lemma~\ref{lemma:projerror}, and assumptions \ref{A3}--\ref{A4}. 
The second term can be bounded by
\begin{align*}
    \eps\|\ptau \hatn w - \dtau\hatn w\|_{L^2} \le C\Dtau
    \|\eps\partial_{\tau\tau}\hat w\|_{L^{\infty}(0,\tau_{max};L^2(0,\ell))},
\end{align*}
and $\|\partial_{\tau\tau}\hat w\|_{L^{\infty}(0,\tau_{max};L^2(0,\ell))}$ can further be estimated by bounds on $\rho$ and $m$ in \ref{A3}--\ref{A4}.
In summary, we thus obtain the following estimate
\begin{align*}
  (i) \le C (h^2 + \Dtau^2) +  C' \H(u_h^n|\hatn u)
\end{align*}
with constants $C,\, C'$ only depending on the bounds in the assumptions.

Using again that \eqref{eq:discr1} and \eqref{eq:pert1} can be tested with any $L^2$-function, we can rewrite 
\begin{align*}
    (ii)=- \la h^n - \hatn h,\dx m_h^n - \dx\hatn m\ra = a \la h^n-\hatn h,\dtau\rho_h^n-\dtau\hatn\rho\ra + \la h^n-\hatn h,\resq^n\ra = (a)+(b).
\end{align*}
Similar to the proof of Lemma \ref{lemma:res1}, we can further bound the second term by
\begin{align*}
    (b) = \la h^n-\hatn h,\resq^n\ra \le \tfrac{1}{2}  \|h^n - \hatn h\|_{L^2}^2 + \tfrac{a^2}{2}\|\dtau\hatn\rho-\partial_{\tau}\hatn\rho\|_{L^2}^2 \le C(\Dtau^2 + h^2),
\end{align*}
where we used that
\begin{align*}
	h^n-\hatn h 
	&\le\tfrac{3}{2} \eps^2\bar w |w^n-\hatn w| +C \, |\rho^n-\hatn\rho|
\le \tfrac{3\eps^2 \bar w}{2 a\ubar\rho}|m^n-\hatn m| + \big(\tfrac{3\eps^2\bar w^2 \bar \rho}{2 \ubar\rho^2} + C\big) |\rho^n-\hatn \rho|\le C'\, h^2
\end{align*}
together with the bounds in the assumptions and Lemma \ref{lemma:projerror}, as well as Lemma \ref{lemma:proj}.
For the first term, we use the following discrete integration-by-parts formula
\begin{align*}
    \dtau u^n v^n = - u^{n-1}\dtau v^n + \tfrac{1}{\Dtau}u^nv^n - \tfrac{1}{\Dtau}u^{n-1}v^{n-1},
\end{align*}
which together with Young's inequality leads to
\begin{align*}
    (a) = \la h^n-&\hatn h,\dtau\rho_h^n-\dtau\hatn\rho\ra
    = -\la \dtau h^n-\dtau\hatn h,\rho_h^{n-1}-\hat\rho_h^{n-1}\ra + \tfrac{1}{\Dtau}\la h-\hat h_h,\rho_h-\hat\rho_h\ra|_{\tau_{n-1}}^{\tau_n}\\
    &\le \tfrac{1}{2}\|\dtau h^n-\dtau\hatn h\|_{L^2}^2 + \tfrac{1}{2}\|\rho_h^{n-1}-\hat\rho_h^{n-1}\|_{L^2}^2 + \tfrac{1}{\Dtau}\la h-\hat h_h,\rho_h-\hat\rho_h\ra|_{\tau_{n-1}}^{\tau_n}.
\end{align*}
The second term can be estimated by $\H(u_h^{n-1},\hat u_h^{n-1})$ and for the first term we use
\begin{align*}
    \|\dtau h^n-\dtau\hatn h\|_{L^2}\le \|\dtau h^n-\ptau h^n\|_{L^2} + \|\ptau h^n-\ptau\hatn h\|_{L^2} + \|\ptau\hatn h - \dtau\hatn h\|_{L^2}.
\end{align*}
We can then further estimate the individual terms by
\begin{align*}
    \|\dtau h^n-\ptau h^n\|_{L^2} &\le C\Dtau\|\partial_{\tau\tau}h\|_{L^\infty(\tau_{n-1},\tau_n;L^2(0,\ell))} \le C' \Dtau,\\
    \|\ptau h^n-\ptau\hatn h\|_{L^2} &\le C'\, h,\\
    \|\ptau\hatn h - \dtau\hatn h\|_{L^2} &\le C \Dtau\|\partial_{\tau\tau}\hatn h\|_{L^\infty(\tau_{n-1},\tau_n;L^2(0,\ell))} \le C' \Dtau,
\end{align*}
with constants that only depend on the bounds in the assumptions. 
For the second inequality in the first and third line, we here used that $\|\partial_{\tau\tau}h\|$ and $\|\partial_{\tau\tau}\hatn h\|$ can be estimated by the bounds for the time derivatives of $\rho$ and $m$ given in \ref{A4}. 
In summary, we thus arrive at 
\begin{align*}
 (ii) \le C (\Dtau^2 + h^2) + \tfrac{1}{\Dtau}\la h-\hat h_h,\rho_h-\hat\rho_h\ra|_{\tau_{n-1}}^{\tau_n} + C'\H(u_h^{n-1}|\hat u_h^{n-1})
\end{align*}
with constants $C,\,C'$ only depending on the bounds in the assumptions.

The remaining term in the estimate \eqref{eq:res2est} can be split into 
 \begin{align*}
 	(iii)&=\big\la \gamma |w^n|w^n-\gamma |\hatn w|\hatn w, \hn m-\hatn m \big\ra
 	=\,\big\la (\gamma |w^n|-\gamma|\hatn w|)(w^n-\hatn w), \hn m-\hatn m \big\ra\\
 	&\qquad +\big\la \gamma |\hatn w| (w^n-\hatn w), \hn m-\hatn m \big\ra
 	+\big\la  \gamma (|w^n|-|\hatn w|) \hatn w, \hn m-\hatn m \big\ra = (A)+(B)+(C).
\end{align*}
By Hölder and Young inequalities, the first term can further be estimated by
\begin{align*}
   (A) =&\ \big\la \gamma (|w^n|-|\hatn w|)(w^n-\hatn w), \hn m-\hatn m \big\ra\\
   \le&\ \bar\gamma \bar a \bar\rho\| (w^n-\hatn w)^2\|_{L^{3/2}} 
   \|w_h^n-\hatn w\|_{L^3} + \tfrac{3}{2} \bar\gamma \bar a \bar w \| (w^n-\hatn w)^2\|_{L^{2}} \|\rho_h^n-\hatn \rho\|_{L^2} \\
   \le&\ \tfrac{2}{3}(\bar \gamma \bar a \bar \rho)^{3/2}\delta^{-3/2} \| w^n-\hatn w\|_{L^{3}}^3
   +\tfrac{1}{3} \delta^3 \| w_h^n-\hatn w\|_{L^{3}}^3
   +C h^2 +C' \H (u_h^n\mid \hatn u)
\end{align*}
for arbitrary $\delta>0$, where we used the projection error bounds of Lemma~\ref{lemma:proj} as well as Lemma~\ref{lemma:relen}. 
Under assumption \ref{A4} the first term in this estimate can be further bounded by $C(\delta)h^3$.
Choosing $\delta$ sufficiently small, we can bound the second term by $\tfrac{1}{4} \D (u_h^n\mid \hatn u)$.
For estimation of the remaining terms we can again use the relative dissipation functional as well as the projection error estimate and Lemma \ref{lemma:relen}, which yields 
\begin{align*}
    (B)+
    (C) &= \big\la \gamma |\hatn w| (w^n-\hatn w), \hn m-\hatn m \big\ra
 	+\big\la \gamma (|w^n|-|\hatn w|) \hatn w, \hn m-\hatn m \big\ra\\
 	&\le \tfrac{1}{\delta'}\bar\gamma^2 \bar w \| w^n-\hatn w\|_{L^{2}}^2
 	+\delta' \|(m_h^n-\hatn m)|\hatn w|^{1/2}\|_{L^{2}}^2\\
    &\le C(\delta') h^2 + 2\delta' \bar a \bar \rho \int_0^\ell a \hatn\rho |w_h^n-\hatn w|^2 |\hatn w|\ dx + \tfrac{9}{2}\delta' \bar a^2\bar w^3 \|\rho_h^n-\hatn \rho\|_{L^2}^2\\
 	&\le C h^2 + \tfrac{1}{4} \D (u_h^n\mid \hatn u) + C' \H (u_h^n\mid \hatn u)
\end{align*}
for $\delta'$ sufficiently small.
In summary, the term $(iii)$ in \eqref{eq:res2est} can be bounded 
by 
\begin{align*}
 	(iii)=&\,\big\la \gamma |w^n|w^n-\gamma |\hatn w|\hatn w, \hn m-\hatn m \big\ra
 	\le C h^2 + C' \H (u_h^n\mid \hatn u) + \tfrac{1}{2} \D (u_h^n\mid \hatn u).
\end{align*} 
Combination of the previous estimates finally yields the assertion of the lemma.
\end{proof}

\noindent
A combination of the bounds in Step 1--3 finally completes the proof of Lemma~\ref{lem:ddH}. 

\subsection{Proof of Theorem~\ref{thm:error_estimate}}\label{ssec:proof}

We can now complete the proof of our main result. 
In view of the error splitting \eqref{eq:split1}--\eqref{eq:split2} and the bounds for the projection error stated in Lemma~\ref{lemma:projerror}, it suffices to consider the discrete error in detail. 
By multiplying \eqref{eq:ddH} with $\Dtau$ and summing over the time steps, we obtain the inequality
\begin{align}  \label{eq:H-estimate}
    \H(u_h^n|\hatn u) \le&\ \H(u_h^0|\hat u_h^0) + \Dtau \sum\nolimits_{k=1}^n\big(C_1\H(u_h^k|\hat u_h^k) + C_2\H(u_h^{k-1}|\hat u_h^{k-1})\big) + (h-\hat h_h)(\rho_h-\hat\rho_h)|_{0}^{\tau_n}\nonumber\\ 
    &\qquad+\Dtau\sum\nolimits_{k=1}^n\big(C_3(\Dtau^2 + h^2) 
    -\tfrac{1}{2}\D(u_h^k|\hat u_h^k) \big).
\end{align}
Using the fact that $\rho_h^0=\hat\rho_h^0$ together with Young's inequality, we can estimate
\begin{align*}
   (h-\hat h_h)(\rho_h-\hat\rho_h)|_{0}^{\tau_n} \le \tfrac{1}{2} c_0^{-1} 
   \|h^n-\hat h_h^n\|^2 + \tfrac{1}{2}\H(u_h^n|\hatn u)
   \le C_4 h^2 + \tfrac{1}{2}\H(u_h^n|\hatn u).
\end{align*}
The last term is moved to the left hand side of \eqref{eq:H-estimate},
and we apply Lemma~\ref{lemma:gronwall} with
\begin{align*}
    a^n=\H(u_h^n\mid\hatn u),\quad c=2\max(C_1,C_2),\quad
    b^n = (2\tau_{\max}C_3+2C_4)(\Dtau^2 + h^2),\quad  
    d^n =  \D(u_h^n|\hat u_h^n).
\end{align*}
Noting that $a^0 = \H(u_h^0\mid\hat u_h^0) = 0$, $b^n\ge 0$ and $n\Dtau\le\tau_{\max}$, we thus obtain
\begin{align*}
	&\H(u_h^n\mid\hatn u) + \Dtau\sum_{k=1}^n \D(u_h^k\mid\hat u_h^k) \le
	C_5(\tau_{\max}) (\Dtau^2+h^2) 
\end{align*}
with constant $C_5$ that only depends on the bounds in \ref{A1}--\ref{A4} but not on $\eps$.
By the equivalence of the norm and the relative energy stated in Lemma~\ref{lemma:relen}, we finally arrive at  
\begin{align}\label{eq:discrerr}
    \|\hn\rho-\hatn\rho\|_{L^2(0,\ell)}^2 + \eps^2\|\hn w-\hatn w\|_{L^2(0,\ell)}^2 +   \sum_{k=1}^n \Dtau \|w_h^k-\hat w_h^k\|_{L^3(0,\ell)}^3
    \leq C\big(
    \Dtau^2 + h^2). 
\end{align}
This estimate initially holds for all $n \le N^*$ given by assumption \ref{A3h}. 
Since the constant~$C$ is independent of $h$, $\Dtau$, and $N^*$, we can now show that this assumption automatically holds for $N^*=N$, if $\Dtau_0$ and $h_0$ are chosen sufficiently small and $h \approx \Dtau$. 

\begin{lemma} \label{lem:drop}
Let \ref{A1}--\ref{A4} hold. Then one can choose $\Dtau_0$ and $h_0$ such that \ref{A3h} holds for all $n \le N^*=N$, $0 < h \le h_0$, $0 < \Dtau \le \Dtau_0$ with $h \approx \Dtau$ and for all $0 \le \eps \le \bar \eps$.
\end{lemma}
\begin{proof}
By an inverse inequality, we can estimate 
\begin{align*}
    \|\hat \rho_h^n - \rho_h^n\|_{L^\infty}^2 \le h^{-1} \|\hat \rho_h^n- \rho_h^n\|_{L^2}^2 \le C (h + \tfrac{\Delta \tau^2}{h}).
\end{align*}
Assuming $\Delta \tau \approx h$ and adding the corresponding projection error, we thus obtain 
\begin{align*}
    \|\rho(\tau^n) - \rho_h^n\|_{L^\infty}^2 \le C \Delta \tau,
\end{align*}
which can be made arbitrarily small. 
Next observe that for $h \approx \Delta\tau$, we get
\begin{align}
    \|\dtau (\hat \rho_h^n - \rho_h^n)\|_{L^2}
    &\le \frac{1}{\Delta \tau} (\|\hat \rho_h^n - \rho_h^n\|_{L^2} + \| \hat \rho_h^{n-1}-\rho_h^{n-1}\|_{L^2}) 
    \le C.
\end{align}
By the triangle inequality and Taylor estimates, we then conclude that 
\begin{align*}
    \|\ptau 
    \hat \rho_h(\tau^n) - \dtau \rho_h^n\|_{L^2} 
    &\le \|\ptau
    \hat \rho_h(\tau^n) - \dtau \hat \rho_h^n\|_{L^2} + \|\dtau (\hat \rho_h^n - \rho_h^n)\|_{L^2} 
    \le C',
\end{align*}
where we used the contraction property for the $L^2$ projection, a Taylor estimate in time, and the previous bound. 
By the commuting diagram property of the projections, we get
\begin{align*}
\|\dx (\hat m_h^n - m_h^n)\|_{L^2} = \|\ptau \hat \rho_h(\tau^n) - \dtau \rho_h^n\|_{L^2} \le C'. 
\end{align*}
From the definition of $w_h^n$ and $\hat w_h^n$ and the uniform bounds for $\rho_h^n$ and $\hat \rho_h^n$, we see that
\begin{align*}
    \|w_h^n - \hat w_h^n\|_{L^\infty}^2 \le C \|\rho_h^n - \hat \rho_h^n\|^2_{L^\infty} + C' \|m_h^n - \hat m_h^n\|^2_{L^\infty} = (i) + (ii).
\end{align*}
Using an inverse inequality, we obtain 
\begin{align*}
    (i) \le C h^{-1}\|\hat \rho_h^n - \rho_h^n\|^2_{L^2} \le C (\Dtau^2/h + h).
\end{align*}
By the multiplicative interpolation inequality, the second term can be estimated by 
\begin{align*}
(ii) &\le  C \|\hat m_h^n - m_h^n\|_{L^2}^2 + C' \|\hat m_h^n - m_h^n\|_{L^2}\|\dx (\hat m_h^n - m_h^n)\|_{L^2}.
\end{align*}
By the definition of $m_h^n$ and $\hat m_h^n$ and the uniform bounds for $\hat \rho_h^n$, $\rho_h^n$ and $\hat w_h^n$, $w_h^n$ for $n \le N^*$, we can further estimate
\begin{align*}
\|\hat m_h^n - m_h^n\|_{L^2} 
&\le C \|\hat w_h^n - w_h^n \|_{L^2} + C'\|\rho_h^n - \hat \rho_h^n\| _{L^2} \\
&\le C'' \Dtau^{-1/3} (\Dtau^2 + h^2)^{1/3} + C''' (\Dtau^2 + h^2)^{1/2},
\end{align*}
where we made use of the estimate \eqref{eq:discrerr}. For $h \approx \Dtau$ and $\Dtau \le \Dtau_0$ sufficiently small, both terms (i) and (ii) can be made as small as desired. As a consequence, $w_h^n$ satisfies the same uniform bounds as $\hat w_h^n$ up to some small perturbation that can be fully controlled by the choice of the mesh size.  By the argument of Lemma~\ref{lem:varh}, the next time step $n=N^*+1$ will then also satisfy $w_h^n \in \AS$, and we can continue applying the argument until $n=N$.
\end{proof}

The above estimates provide the desired bounds for the discrete error in $\rho$ and $w$, which by Lemma~\ref{lem:drop} hold for all $n \le N^*=N$.
Using the relations $\hat m_h^n=a \hat \rho_h^n \hat w_h^n$ and $m_h^n=a \rho_h^n w_h^n$ between the discrete mass fluxes and velocities, one can see that 
\begin{align*}
\eps^2 \|\hat m_h^n - m_h^n\|_{L^2}^2  \le 2 \eps^2 \bar a^2  (\bar \rho^2 \|\hat w_h^n - w_h^n\|_{L^2}^2 + 9/4 \bar w^2 \|\hat \rho_h^n - \rho_h^n\|_{L^2}^2) \le C (\Dtau^2 + h^2).   
\end{align*}
In a similar manner, we can bound 
\begin{align*}
    \|\hat m_h^k - m_h^k\|^3_{L^3} \le 3 \bar a^3 (\bar \rho^3 \|\hat w_h^k - w_h^k\|^3_{L^3} + 27/8 \bar w^3 \|\hat \rho_h^k - \rho_h^k\|^3_{L^3}). 
\end{align*}
By H\"older's inequality and the uniform bounds for $\hat \rho_h^k$ and $\rho_h^k$, we further obtain 
\begin{align*}
    \|\hat \rho_h^k - \rho_h^k\|^3_{L^3} \le 3/2 \bar \rho \|\hat \rho_h^k - \rho_h^k\|^2_{L^2}.
\end{align*}
Together with the previous bounds, we see that
$\sum_{k=1}^n \Dtau \|\hat m_h^n - m_h^n\|^3_{L^3} \le C (\Dtau^2 + h^2)$.
A combination of these estimates for the discrete error and Lemma~\ref{lemma:projerror} for the projection error finally yields the assertion of the theorem.
\qed

\section{Extension to networks}
\label{sec:network}

We now extend our considerations to gas transport problems on pipe networks.
The model equations \eqref{eq:gas1}--\eqref{eq:gas2} are then assumed to hold for each pipe, whereas additional coupling conditions are required at pipe junctions. 
A corresponding variational formulation will be derived and a mixed finite element approximation together with an implicit Euler time stepping is proposed for its numerical solution. 
Since the structure of the  problem and its discretization is very similar to \eqref{eq:var1}--\eqref{eq:var2} and \eqref{eq:discr1}--\eqref{eq:discr2} on a single pipe, our error analysis almost verbatim carries over to networks.

\subsection{Network topology and notation}

The network is described by a directed, connected graph with vertices $v\in\V$ and edges $e\in\E$. Edges correspond to pipes and are identified by intervals $(0,\ell^e)$ for $e\in\E$. Moreover, the edges incident to some vertex $v$ are collected in the set $\E(v)$. We distinguish between boundary vertices $\V_{\partial}=\{v\in\V:|\E(v)|=1\}$ and inner vertices $\V_0=\V\backslash\V_{\partial}$, where $|\E(v)|$ denotes the cardinality of the set $\E(v)$. To each edge $e=(0,\ell^e)=(v_1,v_2)$ we link two numbers $n^e(v_1)=-1,\ n^e(v_2)=1$ to indicate start and end point, and set $n^e(v)=0$ for all $v\in\V\backslash\{v_1,v_2\}$.

We further denote by $L^2(\E)=\{u: u^e\in L^2(e)\ \forall e\in\E\}$ the space of square integrable functions on the network, where $u^e=u|_e$ is the restriction onto the edge $e$. The corresponding scalar product and norm are given by
\begin{align*}
    \la u,v\ra_{\E}:=\sum\nolimits_{e\in\E} \la u^e,v^e\ra_{L^2(e)},\qquad \|u\|_{\E}^2:=(u,u)_{\E}.
\end{align*}
Similarly, we introduce the space $H^1_{pw}(\E)=\{u:u^e\in H^1(e)\}$ of edge-wise $H^1$-functions that are continuous along edges but can be discontinuous at network junctions. We associate the following scalar product and norm
\begin{align*}
    \la u,v\ra_{H^1_{pw}(\E)} := \sum\nolimits_{e\in\E} \la u^e,v^e\ra_{H^1(e)},\qquad \|u\|_{H^1_{pw}(\E)}^2:=\la u,u\ra_{H^1_{pw}(\E)}.
\end{align*}
Other functions spaces can be defined in a similar manner.
\subsection{Gas transport in pipe networks}

We assume that the model equations \eqref{eq:gas1}--\eqref{eq:gas2} are satisfied on each pipe $e\in\E$, i.e.,
\begin{align}
	a^e \partial_\tau \rho^e + \dx m^e &= 0,\label{gasnet:1}\\
	\eps^2 \partial_\tau w^e + \dx h^e &= -\gamma^e|w^e|w^e,\label{gasnet:2}
\end{align}
for $0<x<\ell^e,\ \tau>0$ and all $e\in\E$ with
\begin{align*}
    m^e=a^e\rho^e w^e,\qquad h^e=\tfrac{1}{2}\eps^2(w^e)^2 +P'(\rho^e).
\end{align*}
To guarantee conservation principles at junctions $v \in \V_0$, we impose the coupling conditions
\begin{alignat}{2}
    \sum\nolimits_{e\in\E(v)}m^e(v)n^e(v) &= 0, \qquad &&v\in\V_0,\label{gasnet:3}\\
    h^e(v) &= h^v, &&v\in\V_0,\ e\in\E(v).\label{gasnet:4}
\end{alignat}
For a convenient formulation, the enthalpy $h^v$ has been introduced as an additional degree of freedom for each $v \in \V_0$. As shown in  \cite{Egger18, Reigstad15}, these coupling conditions yield conservation of mass and energy and thus a thermodynamically consistent behavior at pipe junctions.
At the boundary vertices of the network, we again prescribe
\begin{align}
h^e(v)=h_\partial^v,\qquad v\in\V_{\partial},\ e\in\E(v)\label{gasnet:5}. 
\end{align}
A classical solution of \eqref{gasnet:1}--\eqref{gasnet:5} then is a pair of functions 
$$
\rho,w\in C^1([0,\tau_{\max}];L^2(\E)) \cap C^0([0,\tau_{\max}];H^1_{pw}(\E))$$ 
that satisfies the above equations in a pointwise sense. 
In particular, for every point in time, the co-state variable $m=a \rho w$ lies in the space
\begin{align*}
    H(\operatorname{div};\E):=\{u\in H^1_{pw}(\E): \sum\nolimits_{e\in\E}u^e(v)n^e(v)=0\ \forall v\in\V_0\}
\end{align*}
of mass fluxes that are conservative across junctions. 
With a similar reasoning as on a single pipe and the use of the coupling conditions \eqref{gasnet:3}--\eqref{gasnet:4}, one can see that any classical solution satisfies, for all $0 \le \tau \le \tau_{max}$ of interest, the variational identities
\begin{align}
	\la a\partial_{\tau}\rho, q\ra_{\E} + \la \dx m, q\ra_{\E} &= 0 &&\forall q\in L^2(\E), \label{varnet:1}\\
	\la \eps^2\partial_{\tau} w, r\ra_{\E} - \la h, \dx r\ra_{\E} + \la \gamma |w| w, r\ra_{\E}&= - \sum\nolimits_{v\in\V_{\partial}} h_\partial^v\, r^e(v)n^e(v) &&\forall r\in H(\operatorname{div};\E). \label{varnet:2}
\end{align}

\begin{remark}
Let us note that the coupling condition \eqref{gasnet:3} on the mass flux is strongly enforced in the space $H(\operatorname{div};\E)$, whereas the continuity condition \eqref{gasnet:4} is included in the variational formulation. More precisely, when applying integration by parts on the second term in equation \eqref{gasnet:2} the boundary contributions at inner vertices vanish, i.e.,
\begin{align*}
    \la \dx h,r\ra_{\E} = -\la h,\dx r\ra_{\E} + \sum\nolimits_{v\in\V}\sum\nolimits_{e\in\E(v)} h^e(v)\, r^e(v)n^e(v),
\end{align*}
and since $h$ is assumed to be continuous along junctions and $r\in H(\operatorname{div};\E)$ we see that
\begin{align*}
    \sum\nolimits_{v\in\V_0}\sum\nolimits_{e\in\E(v)} h^e(v)r^e(v)n^e(v) = \sum\nolimits_{v\in\V_0}h^v\,\sum\nolimits_{e\in\E(v)} r^e(v)n^e(v) = 0,
\end{align*}
and only the contributions at the boundary vertices remain. 
\end{remark}
The total energy contained in the network is now simply defined by accumulation of the contributions of the individual pipes, i.e.
\begin{align*}
   \H(\rho,w)=\sum\nolimits_{e\in\E} \int_0^{\ell^e} a^e (\tfrac{1}{2}\eps^2 \rho^e (w^e)^2 +P(\rho^e))\ dx. 
\end{align*}
Similar as on a single pipe, we can again deduce an energy-dissipation law
\begin{align*}
    \ptau \H(\rho,w) + \D(\rho,w) = - \sum\nolimits_{v\in\V_{\partial}} h_\partial^v\, m^e(v)n^e(v)
\end{align*} 
with dissipation functional $\D(\rho,w) = \sum\nolimits_{e\in\E}\int_0^{\ell^e} a^e \gamma^e\rho^e|w^e|^3 \geq 0$, 
which again follows directly from the particular form of the variational formulation.
Based on relative energy estimates, the stability of solutions to \eqref{gasnet:1}--\eqref{gasnet:5} with respect to perturbations in the initial conditions and the problem parameters has been analysed in \cite{EggerGiesselmann}.  
Here we use a similar reasoning to extend our discretization scheme and error estimates to gas networks.

\subsection{A Galerkin scheme on networks}

We approximate the density $\rho$
by piecewise constant functions 
over the grid $\T_h =\{[x_{i-1}^e,x_i^e]: x_0^e = 0,\ x_{M^e}^e=\ell^e,\ x_i^e-x_{i-1}^e = h^e,\ e\in\E\}$ and set $h=\max h^e$.
The mass flux $m$ is approximated by piecewise linear functions over the mesh $\T_h$, which are continuous on every pipe $e$ and satisfy the balance condition \eqref{gasnet:3}.  
The spatial approximation spaces are then given by
\begin{align*}
    Q_h:=\mathcal{P}_0(\T_h),\qquad R_h:=\mathcal{P}_1(\T_h)\cap H(\operatorname{div};\E).
\end{align*}
By $\Pi_h:L^2(\E)\rightarrow Q_h,\ I_h:H(\operatorname{div};\E)\rightarrow R_h$ we denote the canonical extensions of the locally defined projection and interpolation operators to the network setting. 
Based on the variational formulation \eqref{varnet:1}--\eqref{varnet:2} we then propose the following method.
\begin{problem}[Discretization scheme for gas networks]\label{prob:fds-net} $ $\\
Set $\rho_h^0 = \Pi_h\rho(0),\ m_h^0 = I_h m(0)$, and for $1 \le n \le N$ find $\rho_h^n\in Q_h,\ m_h^n\in R_h$ such that
\begin{align}
	\la a \dtau \rho_h^n, q_h\ra_\E + \la \dx m_h^n, q_h\ra_\E =&\ 0  \label{eq:discrnet1} \\
	\la \eps^2 \dtau w_h^n, r_h\ra_\E - \la h_h^n, \dx r_h\ra_\E 
	+ \la \gamma |w_h^n|w_h^n, r_h\ra_\E =&\ -\!\! \sum\nolimits_{v\in\V_{\partial}} \sum\nolimits_{e\in\E(v)} h_\partial^v\, r_h^e(v)n^e(v) \label{eq:discrnet2}
\end{align}
for all $q_h \in Q_h$ and $r_h \in R_h$. 
We again use $w_h^n
=\tfrac{m_h^n}{a\rho_h^n}$ and $ h_h^n
=\tfrac{\eps^2(m_h^n)^2}{2a^2(\rho_h^n)^2}+P'(\rho_h^n)$ to abbreviate the discrete velocity and enthalpy variables.
\end{problem}

\begin{remark}
With the very same arguments as on a single pipe, one can show that solutions of Problem~\ref{prob:fds-net} satisfy the discrete energy inequality
\begin{align*}
    \dtau \H(\rho_h^n,w_h^n) + \D(\rho_h^n,w_h^n) \leq - \sum\nolimits_{v\in\V_{\partial}} \sum\nolimits_{e\in\E(v)} h_\partial^v m_h^e(v)n^e(v)
\end{align*} 
with dissipation functional $\D(\rho,w) = \sum\nolimits_{e\in\E}\int_0^{\ell^e} a^e \gamma^e\rho^e|w^e|^3 \geq 0$.
\end{remark}

\subsection{Error analysis}

Since the structure of the variational problem is exactly the same as for a single pipe, the analysis of the previous sections carries over verbatim by simply summing over all pipes.
As an immediate consequence, we obtain the following result.
\begin{theorem}\label{thm:error_estimate_net}
Let the assumptions of Theorem~\ref{thm:error_estimate} hold for the network setting. Then
\begin{align*}
	\|\rho(\tau^n)-\rho_h^n\|_{L^2(\E)}^2 + \eps^2\|m(\tau^n)-m_h^n\|_{L^2(\E)}^2 +   \sum_{k=1}^n \Dtau \|m(\tau^k)-m_h^k\|_{L^3(\E)}^3
	\leq C\big(\Dtau^2 + h^2\big), 
\end{align*}
with $C$ depending only on the bounds in the assumptions, but independent of $\eps$.
\end{theorem}


\section{Numerical experiments}
\label{sec:numerics}

For illustration of our theoretical results, we now report about some numerical tests. 
In the first example, we consider the $\eps$-robustness of the convergence estimates for a single pipe, and in the second example, we briefly address the extension to pipe networks. 

\subsection{Parameter robust convergence}

We consider the flow through pipes of different lengths $L \approx \eps^{-2}$ with the diameter and friction coefficient kept fixed. 
By the rescaling outlined in Appendix~\ref{app:rescale}, we can transform the equations into the system \eqref{eq:gas1}--\eqref{eq:gas2} for a rescaled pipe of length $\ell=1$, with uniform cross section and friction coefficient, but with different scaling parameters $\eps$. 
For ease of presentation, we set $a=1$ and  $\gamma=1$, and as a pressure law, we choose $p(\rho)=c^2 \rho$, with speed of sound rescaled to $c=1$. 
As boundary conditions, we choose
\begin{align*}
h_\partial^0(\tau)=0.2\sin(\pi \tau)^3+1, \qquad h_\partial^\ell(\tau)=0.1\sin(\pi+\pi \tau)^3+1
\end{align*}
over a time horizon of $\tau_{max}=1$.
The initial conditions are determined by solving the stationary problem for the boundary conditions at time $\tau=0$. 
Let us note that for $\eps \ll 1$, we have $h =\eps^2 \frac{w^2}{2} + P'(\rho) \approx P'(\rho)$; hence setting the enthalpy is more or less equivalent to prescribing the density or the pressure, respectively. 
In Table \ref{tab:1} we display the errors and convergence rates in density and mass flux for different choices of the scaling parameter $\eps$.
Since the exact solution is unknown, the numerical errors are computed as
\begin{align}
    \text{err}_h(u) = \max_{n=1,..,N} \|u_h^n-u_{h/2}^n\|_{L^2(0,1)},\label{err:0}
\end{align}
where $u = \rho$ or $u=m$ and with $u_{h/2}^n$ denoting the solution on a finer mesh with $h=h/2$ and $\Dtau=\Dtau/2$ at the same point $\tau_n=n \Dtau$ in time.
\begin{table}[ht]
\centering
\footnotesize
\begin{tabular}{c|c|c|c|c|c|c|c}
\multirow{4}{*}{$\eps=1$} 
& err$_h(\rho)$ & 1.28e-2 & 7.58e-3 & 4.21e-3 & 2.24e-3 & 1.16e-3 & 5.89e-4\\
& rate & --- & 0.76 & 0.85 & 0.91 & 0.95 & 0.97\\
\cline{2-8}
& err$_h(m)$ & 1.17e-2 & 7.19e-3 & 4.06e-3 & 2.19e-3 & 1.15e-3 & 5.92e-4\\
& rate & --- & 0.71 & 0.83 & 0.89 & 0.93 & 0.96\\
\hline\hline
\multirow{4}{*}{$\eps=0.1$} 
& err$_h(\rho)$ & 4.99e-3 & 2.49e-3 & 1.25e-3 & 6.23e-4 & 3.12e-4 & 1.56e-4\\
& rate & --- & 1.00 & 1.00 & 1.00 & 1.00 & 1.00\\
\cline{2-8}
& err$_h(m)$ & 9.61e-3 & 5.47e-3 & 2.92e-3 & 1.52e-3 & 7.79e-4 & 3.93e-4\\
& rate & --- & 0.81 & 0.90 & 0.94 & 0.97 & 0.98\\
\hline\hline
\multirow{4}{*}{$\eps=0.01$} 
& err$_h(\rho)$ & 4.98e-3 & 2.49e-3 & 1.24e-3 & 6.22e-4 & 3.11e-4 & 1.55e-4\\
& rate & --- & 1.00 & 1.00 & 1.00 & 1.00 & 1.00\\
\cline{2-8}
& err$_h(m)$ & 4.10e-3 & 2.09e-3 & 1.06e-3 & 5.32e-4 & 2.95e-4 & 1.56e-4\\
& rate & --- & 0.97 & 0.99 & 0.99 & 0.85 & 0.92\\
\hline\hline
\multirow{4}{*}{$\eps=0.001$} 
& err$_h(\rho)$ & 4.98e-3 & 2.49e-3 & 1.24e-3 & 6.22e-4 & 3.11e-4 & 1.55e-4\\
& rate & --- & 1.00 & 1.00 & 1.00 & 1.00 & 1.00\\
\cline{2-8}
& err$_h(m)$ & 4.10e-3 & 2.09e-3 & 1.06e-3 & 5.31e-4 & 2.66e-4 & 1.33e-4\\
& rate & --- & 0.97 & 0.99 & 0.99 & 1.00 & 1.00\\
\hline\hline
\multirow{4}{*}{$\eps=0$} 
& err$_h(\rho)$ & 4.98e-3 & 2.49e-3 & 1.24e-3 & 6.22e-4 & 3.11e-4 & 1.55e-4\\
& rate & --- & 1.00 & 1.00 & 1.00 & 1.00 & 1.00\\
\cline{2-8}
& err$_h(m)$ & 4.10e-3 & 2.09e-3 & 1.06e-3 & 5.31e-4 & 2.66e-4 & 1.33e-4\\
& rate & --- & 0.97 & 0.99 & 0.99 & 1.00 & 1.00\\
\end{tabular}
\vskip1em
\caption{Error and convergence rates for different values of $\eps$. Space and time discretization with $h=\frac{1}{16} \cdot 2^{-r}$ and $\Delta\tau = \tfrac{1}{2} h$ in refinement $r=0,\ldots,5$.}
\label{tab:1}
\end{table}
As predicted by our theoretical results, we observe linear convergence for both density and mass flux uniform for all parameters $\eps$ and, in particular, also in the parabolic limit $\eps=0$. 
Further note that the errors, and actually also the solutions, are very similar for all values of $\eps \le 0.01$, which clearly indicates the asymptotic convergence of solutions with $\eps \searrow 0$, which was proven in \cite{EggerGiesselmann} for the continuous problem.

\subsection{A simple gas network}
\begin{figure}[ht!]
\centering
\begin{tiny}
\begin{tikzpicture}[scale=2]
\node (A) at (0,0) [circle,draw,thick] {$v_1$};
\node (B) at (1,0) [circle,draw,thick] {$v_2$};
\node (D) at (2,0) [circle,draw,thick] {$v_2'$};
\node (G) at (2.75,-0.5) [circle,draw,thick] {$v_2''$};
\node (E) at (2.75,0.5) [circle,draw,thick] {$v_3$};
\node (F) at (2.75,1.2) [circle,draw,thick] {$v_4$};
\node (H) at (2.75,-1.2) [circle,draw,thick] {$v_5$};
\node (I) at (3.5,0) [circle,draw,thick] {$v_6$};
\node (J) at (4.5,0) [circle,draw,thick] {$v_6'$};
\node (K) at (5.25,0.5) [circle,draw,thick] {$v_7$};
\node (L) at (5.25,-0.5) [circle,draw,thick] {$v_8$};
\draw[->, very thick] (A) to node[above] {$e_1$} (B);
\draw[->, very thick] (D) to node[above left=-0.1cm] {$e_2$} (E);
\draw[->, very thick] (H) to node[right] {$e_3$} (G);
\draw[->, very thick] (E) to node[right] {$e_4$} (F);
\draw[->, very thick] (E) to node[above right=-0.1cm] {$e_5$} (I);
\draw[->, very thick] (G) to node[above left=-0.1cm] {$e_6$} (I);
\draw[->, very thick] (J) to node[above] {$e_7$} (K);
\draw[->, very thick] (J) to node[above] {$e_8$} (L);
\draw[->, very thick,densely dashed,black!60] (B) to node[above] {$e_{cs1}$} (D);
\draw[->, very thick,densely dashed,black!60] (I) to node[above] {$e_{cs2}$} (J);
\draw[->, very thick,densely dashed,black!60] (D) to node[above right=-0.1cm] {$e_{vlv}$} (G);
\end{tikzpicture}
\end{tiny}
\caption{GasLib-11 network.}
\label{fig:gaslib11}
\end{figure}
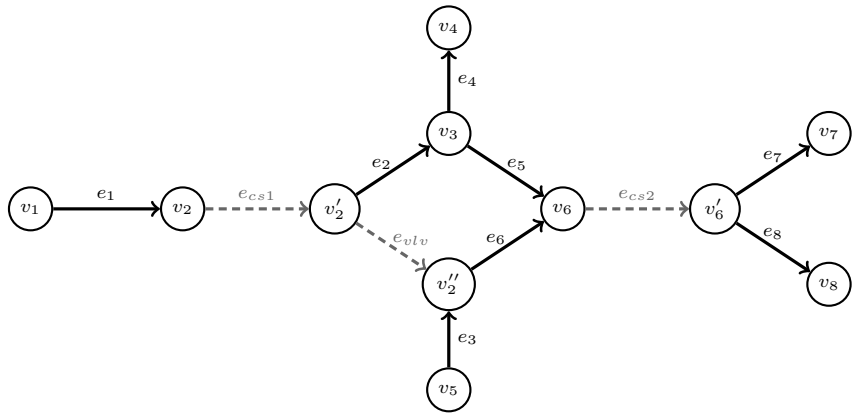

As a second example, we consider the GasLib-11 example from the GasLib library \cite{Gaslib}; see Figure~\ref{fig:gaslib11} for a sketch of the network topology.  
This network consists of $8$ pipes and $11$ vertices, $4$ of which are boundary vertices in the sense of the previous section. The $2$ compressor stations and the valve  are treated here in \emph{by-pass mode}, i.e., as additional pipes of length $0$. The network thus degenerates to a reduced network, where the vertices $v_2,\,v_2',\ v_2''$ as well as $v_6,\,v_6'$ are identified by single vertices $v_2$ and $v_6$ with ingoing edges $e_1,\, e_3$ and $e_5,\, e_6$ as well as outgoing edges $e_2,\,e_6$ and $e_7,\,e_8$, respectively.
All pipes are assumed to have rescaled length $\ell=1$ with cross-sectional area $a=1$ and friction coefficient $\gamma=1$. The pressure law is again given by $p(\rho)=c^2\rho$ with rescaled speed of sound $c=1$. The boundary conditions for the enthalpy are chosen as
\begin{align*}
    h_{\partial}^{v_1}(\tau) = 0.2 \sin(\pi\tau)^3+1,\ \ h_{\partial}^{v_5}(\tau) = 0.3 \sin(\pi\tau)^3+1,\ \ h_{\partial}^{v_4}(\tau)=h_{\partial}^{v_7}(\tau)=h_{\partial}^{v_8}(\tau)=1
\end{align*}
over a time horizon of $\tau_{max}=1$, and the initial condition is again given by the stationary state corresponding to the boundary conditions at $\tau=0$. 
\begin{table}[ht]
\centering
\footnotesize
\begin{tabular}{c|c|c|c|c|c|c|c}
\multirow{4}{*}{$\eps=1$} 
& err$_h(\rho)$ & 2.01e-2 & 1.31e-2 & 8.07e-3 & 4.64e-3 & 2.54e-3 & 1.34e-3\\
& rate & - & 0.61 & 0.70 & 0.80 & 0.87 & 0.92\\
\cline{2-8}
& err$_h(m)$ & 1.70e-2 & 1.11e-2 & 6.72e-3 & 3.85e-3 & 2.11e-3 & 1.11e-3\\
& rate & - & 0.61 & 0.72 & 0.80 & 0.87 & 0.92\\
\hline\hline
\multirow{4}{*}{$\eps=10^{-1}$} 
& err$_h(\rho)$ & 6.01e-3 & 3.04e-3 & 1.57e-3 & 8.23e-4 & 4.32e-4 & 2.24e-4\\
& rate & - & 0.98 & 0.96 & 0.93 & 0.93 & 0.95\\
\cline{2-8}
& err$_h(m)$ & 2.90e-2 & 1.74e-2 & 9.83e-3 & 5.34e-3 & 2.83e-3 & 1.47e-3\\
& rate & - & 0.74 & 0.82 & 0.88 & 0.92 & 0.94\\
\hline\hline
\multirow{4}{*}{$\eps=10^{-2}$} 
& err$_h(\rho)$ & 6.04e-3 & 3.03e-3 & 1.52e-3 & 7.60e-4 & 3.80e-4 & 1.90e-4\\
& rate & - & 0.99 & 1.00 & 1.00 & 1.00 & 1.00\\
\cline{2-8}
& err$_h(m)$ & 2.04e-2 & 1.29e-2 & 7.71e-3 & 4.13e-3 & 2.20e-3 & 1.14e-3\\
& rate & - & 0.66 & 0.74 & 0.90 & 0.91 & 0.95\\
\hline\hline
\multirow{4}{*}{$\eps=10^{-3}$} 
& err$_h(\rho)$ & 6.04e-3 & 3.03e-3 & 1.52e-3 & 7.61e-4 & 3.80e-4 & 1.90e-04\\
& rate & - & 0.99 & 1.00 & 1.00 & 1.00 & 1.00\\
\cline{2-8}
& err$_h(m)$ & 2.11e-2 & 1.31e-2 & 7.30e-3 & 3.95e-3 & 2.05e-3 & 1.05e-03\\
& rate & - & 0.69 & 0.84 & 0.89 & 0.94 & 0.97\\
\hline\hline
\multirow{4}{*}{$\eps=0$} 
& err$_h(\rho)$ & 6.04e-3 & 3.03e-3 & 1.52e-3 & 7.61e-4 & 3.81e-4 & 1.90e-04\\
& rate & - & 0.99 & 1.00 & 1.00 & 1.00 & 1.00\\
\cline{2-8}
& err$_h(m)$ & 2.12e-2 & 1.31e-2 & 7.30e-3 & 3.95e-3 & 2.05e-3 & 1.05e-3\\
\cline{2-8}
& rate & - & 0.69 & 0.84 & 0.89 & 0.94 & 0.97\\
\end{tabular}
\vskip1em
\caption{Error and convergence rates for different values of $\eps$. Space and time discretization with $h=\frac{1}{16} \cdot 2^{-r}$ and $\Delta\tau = \frac{1}{2} h$ in refinement $r=0,\ldots,5$.}
\label{tab:2}
\end{table}
In Table~\ref{tab:2} errors in density and mass flux, computed as in \eqref{err:0}, and convergence rates are presented. Again, we observe linear convergence uniform for all choices of the parameter $\eps$. 
In the light of our theoretical results, the numerical results for the network are expected and also observed to be quasi identical to those for a single pipe; compare Table~\ref{tab:1} and \ref{tab:2}.

According to the specifications in \cite{Gaslib}, the real length of the pipes in the GasLib-11 network is about $55$km and their diameter about $0.5$m. Interesting time scales are in the range of hours and days, which corresponds to a scaling parameter of $\eps\approx 0.01-0.001$ in the parabolic limit regime; see \cite{Brouwer11} for details. 
As observed in \cite{Osiadacz84} and illustrated in our numerical tests, one can set $\eps=0$ in that case and directly use the parabolic limit problem for simulations with practically the same outcome. 
See \cite{Bamberger79,Burlacu19,SchoebelKroehn20} for alternative results for the parabolic limit problem obtained by a different discretization strategy.

\section*{Acknowledgement}
The authors are grateful for financial support by the German Science Foundation (DFG) via grant TRR~154 (\emph{Mathematical modelling, simulation and optimization using the example of gas networks}), projects~C04 and C05.


\begin{thebibliography}{10}


\bibitem{Brouwer11}
J.~Brouwer, I.~Gasser, and M.~Herty.
\newblock Gas pipeline models revisited: model hierarchies, nonisothermal models, and simulations of networks.
\newblock {\em Multiscale Model. Simul.}, 9:601--623, 2011.

\bibitem{Egger18}
H.~Egger.
\newblock A robust conservative mixed finite element method for isentropic
  compressible flow on pipe networks.
\newblock {\em SIAM J. Sci. Comput.}, 40:A108--A129, 2018.

\bibitem{Reigstad15}
G.~A. Reigstad.
\newblock Existence and uniqueness of solutions to the generalized {R}iemann
  problem for isentropic flow.
\newblock {\em SIAM J. Appl. Math.}, 75:679--702, 2015.

\bibitem{EggerGiesselmann}
H.~Egger and J.~Giesselmann.
\newblock Stability and asymptotic analysis for instationary gas transport via
  relative energy estimates.
\newblock {\em arXiv:2012.14135}, 2020.

\bibitem{GiesselmannLattanzioTzavaras_17}
J.~{Giesselmann}, C.~{Lattanzio}, and A.~E. {Tzavaras}.
\newblock {Relative energy for the {K}orteweg theory and related {H}amiltonian
  flows in gas dynamics}.
\newblock {\em Arch. Ration. Mech. Anal.}, 223:1427--1484, 2017.

\bibitem{Egger19}
H.~Egger.
\newblock Structure preserving approximation of dissipative evolution problems.
\newblock {\em Numer. Math.}, 143:85--106, 2019.

\bibitem{Geveci88}
T.~Geveci.
\newblock On the application of mixed finite element methods to the wave
  equations.
\newblock {\em RAIRO Mod\'{e}l. Math. Anal. Num\'{e}r.}, 22:243--250, 1988.

\bibitem{Joly03}
P.~Joly.
\newblock Variational methods for time-dependent wave propagation problems.
\newblock In {\em Topics in computational wave propagation}, volume~31 of {\em Lect. Notes Comput. Sci. Eng.}, pp.~201--264. Springer, 
  2003.

\bibitem{Cardoso19}
F.~L. Cardoso-Ribeiro, D.~Matignon, and L.~Lef{\`e}vre.
\newblock A partitioned finite element method for power-preserving
  discretization of open systems of conservation laws.
\newblock {\em arXiv:1906.05965}, 2019.

\bibitem{Liljegren20}
B.~Liljegren-Sailer and N.~Marheineke.
\newblock On port-{H}amiltonian approximation of a nonlinear flow problem on
  networks.
\newblock {\em arXiv:2009.11216}, 2020.

\bibitem{Burlacu19}
R.~Burlacu, H.~Egger, M.~Gro{\ss}, A.~Martin, M.~E. Pfetsch, L.~Schewe,
  M.~Sirvent, and M.~Skutella.
\newblock Maximizing the storage capacity of gas networks: a global {MINLP}
  approach.
\newblock {\em Optim. Eng.}, 20:543--573, 2019.

\bibitem{Bamberger79}
A.~Bamberger, M.~Sorine, and J.~P. Yvon.
\newblock Analyse et contr\^{o}le d'un r\'{e}seau de transport de gaz.
\newblock In {\em Computing methods in applied sciences and engineering
  (Proc. Third Internat. Sympos., Versailles, 1977), {II}}, volume~91
  of {\em Lecture Notes in Phys.}, pp.~347--359. Springer, Berlin-New York, 1979.

\bibitem{SchoebelKroehn20}
L.~Sch\"obel-Kr\"ohn.
\newblock {\em Analysis and numerical approximation of nonlinear evolution
  equations on network structures}.
\newblock Dr. Hut-Verlag, M\"unchen, 2020.

\bibitem{hertyEMS}
A.~Bressan, S.~\v{C}ani\'c, M.~Garavello, M.~Herty, and B.~Piccoli.
\newblock Flows on networks: recent results and perspectives.
\newblock {\em EMS Surv. Math. Sci.}, 1:47--111, 2014.

\bibitem{Jungel16}
A.~J{\"u}ngel.
\newblock {\em Entropy methods for diffusive partial differential equations}.
\newblock Springer, 2016.

\bibitem{Dafermos05}
C.~M. Dafermos.
\newblock {\em Hyperbolic conservation laws in continuum physics}.
\newblock Springer, 2005.

\bibitem{Feireisl18}
E.~Feireisl, M.~Lukacova-Medvidova, S.~Necasova, N.~Antonin, and B.~She.
\newblock Asymptotic preserving error estimates for numerical solutions of
  compressible {N}avier--{S}tokes equations in the low {M}ach number regime.
\newblock {\em Multiscale Model. Simul.}, 16:150--183, 2018.

\bibitem{Gallouet16}
T.~Gallou{\"e}t, R.~Herbin, D.~Maltese, and A.~Novotny.
\newblock Error estimates for a numerical approximation to the compressible
  barotropic {N}avier--{S}tokes equations.
\newblock {\em IMA J. Numer. Anal.}, 36:543--592, 2016.

\bibitem{Kwon20}
Y.-S. Kwon and A.~Novotny.
\newblock Consistency, convergence and error estimates for a mixed finite
  element/finite volume scheme to compressible {N}avier-{S}tokes equations with
  general inflow/outflow boundary data.
\newblock {\em arXiv:2005.00799}, 2020.

\bibitem{Berthon_Bessemoulin_Mathis_2017}
C.~Berthon, M.~Bessemoulin-Chatard, and H.~Mathis.
\newblock Numerical convergence rate for a diffusive limit of hyperbolic
  systems: \(p\)-system with damping.
\newblock {\em SMAI J. Comput. Math.}, 2:99--119, 2016.

\bibitem{Yee_1981}
H.~C. Yee.
\newblock Numerical approximation of boundary conditions with applications to
  inviscid equations of gas dynamics.
\newblock Technical Report TM-18265, NASA, 1981.

\bibitem{BrennerScott_2008}
S.~C. Brenner and L.~R. Scott.
\newblock {\em The mathematical theory of finite element methods}, volume~15 of
  {\em Texts in Applied Mathematics}.
\newblock Springer, New York, 2008.

\bibitem{Gaslib}
M.~Schmidt, D.~A{\ss}mann, R.~Burlacu, J.~Humpola, I.~Joormann, N.~Kanelakis,
  T.~Koch, D.~Oucherif, M.~E. Pfetsch, L.~Schewe, R.~Schwarz, and M.~Sirvent.
\newblock {G}as{L}ib -- {A} {L}ibrary of {G}as {N}etwork {I}nstances.
\newblock {\em Data}, 2, 2017.

\bibitem{Osiadacz84}
A.~Osiadacz.
\newblock Simulation of transient gas flows in networks.
\newblock {\em Int. J. Numer. Meth. Fluids}, 4:13--24, 1984.

\end{thebibliography}


\newpage

\appendix

\renewcommand{\thetheorem}{A\arabic{theorem}}

\section{Transformation and rescaling of model equations}\label{app:rescale}

Consider the one dimensional barotropic Euler equations with friction \begin{align}
    a\dt\rho + \dx m =& 0,\label{oeq:gas1}\\
    \dt m + \dx(\tfrac{m^2}{a\rho} + a p(\rho)) =& -\tfrac{\lambda}{2d}\tfrac{|m|m}{a\rho} v\label{oeq:gas2}
\end{align}
with gas density $\rho$, mass flux $m=a\rho v$, flow velocity $v$, pipe diameter $d$ and cross-sectional area $a$, and friction coefficient $\lambda$.
Using the product rule of differentiation in \eqref{oeq:gas2} together with \eqref{oeq:gas1}, one can see that
\begin{align*}
    \dt v 
    &= -\tfrac{1}{a\rho}\big(av\dt\rho + v\dx(a\rho v) +a\rho v\dx v + a\dx p(\rho) + \tfrac{\lambda}{2d}|v|a\rho v\big) \\
    &=-\tfrac{1}{2}\dx v^2 - 
    \tfrac{1}{\rho} \dx p(\rho) 
    - \tfrac{\lambda}{2d}|v|v.
\end{align*}
We further introduce the pressure potential 
$P(\rho) = \rho\int_1^\rho \frac{p(r)}{r^2}\ dr$
and observe that 
\begin{align*}
    \dx P'(\rho) = \dx\big(\tfrac{p(\rho)}{\rho} + \int_1^\rho \tfrac{p(r)}{r^2}\ dr\big)
    = -\tfrac{p(\rho)}{\rho^2}\dx\rho + \tfrac{\dx p(\rho)}{\rho} + \tfrac{p(\rho)}{\rho^2}\dx\rho = \tfrac{1}{\rho} \dx p(\rho),
\end{align*}
which allows to rewrite the evolution equation for the velocity compactly as
\begin{align} \label{eq:v}
\dt v + \dx (\tfrac{v^2}{2} + P'(\rho)) = - \tfrac{\lambda}{2d} |v| v. 
\end{align}

We then employ two rescalings of the model equations:
In a first step, we replace
\begin{align} \label{eq:rescale_initial}
x \to x/\eps^2, \qquad t \to t/\eps^2,
\end{align}
which resembles the situation of long pipes and time scales. 
After division by $\eps^2$, this leads to a friction term with parameter $\lambda/\eps^2 \to \lambda=O(1/\eps^2)$, characterizing the  \emph{large friction regime}.
In a second step, we rescale  \eqref{oeq:gas1} and \eqref{oeq:gas2}, now with parameter $\lambda=O(1/\eps^2)$, by
\begin{align} \label{eq:rescale}
    t=\tfrac{1}{\eps}\tau,\qquad 
    v=\eps w,\qquad
    \lambda = \tfrac{2 d}{\eps^2} \gamma, 
\end{align}
which corresponds to the long time, small velocity  and low Mach setting of relevance in the large friction case $\eps \ll 1$ typical for the gas transport in long pipelines. 

A direct application of the these rescalings to \eqref{oeq:gas1}--\eqref{oeq:gas2} leads to the system \eqref{eq:cgas_1}--\eqref{eq:cgas_2} considered in the introduction,
and together with the above transformation of the momentum equation, one obtains the system \eqref{eq:gas1}--\eqref{eq:gas2} considered in Section~\ref{sec:single}.


\newpage 

\renewcommand{\thetheorem}{B\arabic{theorem}}

\section{Discrete Gronwall lemma}  \label{app:gronwall}

For the proof of the discrete stability estimate, we employ the following technical result.
\begin{lemma}[Discrete Gronwall]\label{lemma:gronwall}
Let $a^n$, $b^n$, $d^n \ge 0$ for $n=0,\ldots,N$, and further  $N \Dtau  = \tau_{\max}$  and $c>0$ with $c\Dtau<1$ be given, such that
\begin{align}\label{eq:gron1}
    a^n +  \sum_{k=1}^n \Dtau d^k \le a^0 +  b^n + c \sum_{k=1}^n \Dtau (a^k + a^{k-1}).
\end{align}
Then, it holds
\begin{align*}
    a^n + \sum_{k=1}^n \Dtau d^k \le a^0 + b^n + c\Dtau e^{\frac{2nc\Dtau}{1-c\Dtau}}\big(a^0 + \sum_{k=1}^n e^{\frac{(1-2k)c\Dtau}{1-c\Dtau}}(2a_0 + b^k + b^{k-1})\big).
\end{align*}
\end{lemma}
\begin{proof}
We set 
\begin{align*}
    s^n:=\sum\nolimits_{k=1}^n (a^k+a^{k-1}),\quad s^0=a^0\quad \text{and}\quad w:=\tfrac{1-c\Dtau}{1+c\Dtau}.
\end{align*}
Then \eqref{eq:gron1} can be written as
\begin{align}\label{eq:gron3}
    a^n - c\Dtau s^n  \le a^0 + b^n  \sum\nolimits_{i=k}^n \Dtau d^k.
\end{align}
Now, define $\tilde a^n:= w^ns^n$. Then
\begin{align*}
    \tilde a^n-\tilde a^{n-1} 
    &= w^n s^n-w^{n-1}s^{n-1} = w^{n-1}(ws^n-s^{n-1})\\
    &= w^{n-1}(1+c\Dtau)^{-1}(s^n-c\Dtau s^n-s^{n-1}-c\Dtau s^{n-1})\\
    &= w^{n-1}(1+c\Dtau)^{-1}\big( (a^n -c\Dtau s^n)+(a^{n-1}-c\Dtau s^{n-1}) \big)\\   
    &\le w^{n-1}(1+c\Dtau)^{-1}(a^0 + b^n - \sum\nolimits_{k=1}^n \Dtau d^k+a^0 + b^{n-1} -  \sum\nolimits_{k=1}^{n-1} \Dtau d^k\big),
\end{align*}
where we used \eqref{eq:gron3} in the last estimate. 
Summing up over $n$ yields 
\begin{align*}
    \tilde a^n\le \tilde a^0 + \sum\nolimits_{k=1}^n w^{k-1}(1+c\Dtau)^{-1}\big( 2a^0 + b^k + b^{k-1} - \sum\nolimits_{j=1}^{k} \Dtau d^j - \sum\nolimits_{j=1}^{k-1} \Dtau d^j\big).
\end{align*}
From the definitions $\tilde a^0 = a^0$ and $\tilde a^n = w^ns^n$ as well as $d^j\ge 0$, we then deduce that
\begin{align*}
    s^n \le w^{-n}a^0 + w^{-n}\sum\nolimits_{k=1}^n w^{k-1}(1+c\Dtau)^{-1}\big( 2a^0 + b^k + b^{k-1}\big).
\end{align*}
For the terms on the right hand side, we further use 
\begin{align*}
    w^{-n}w^{k-1}(1+c\Dtau)^{-1} = \big(\tfrac{1-c\Dtau}{1+ c\Dtau}\big)^{k-n-1}(1+c\Dtau)^{-1} \le e^{\frac{(2(n-k)+1)c\Dtau}{1-c\Dtau}},
\end{align*}
which then leads to 
\begin{align*}
    a^n + \Dtau\sum_{k=1}^n d^k \le a^0 + b^n + c\Dtau e^{\frac{2nc\Dtau}{1-c\Dtau}}\big(a^0 + \sum_{k=1}^n e^{\frac{(1-2k)c\Dtau}{1-c\Dtau}}(2a_0 + b^k + b^{k-1})\big).
\end{align*}
With $n\Dtau\le N \Dtau = \tau_{\max}$, we finally obtain the claim of the lemma. 
\end{proof}

\end{document}